\documentclass[pdflatex,sn-mathphys-num]{sn-jnl}
\usepackage{amsmath,amssymb,amsthm}
\usepackage{mathabx}
\usepackage{graphicx} 
\usepackage{subcaption}
\usepackage{enumerate}
\usepackage{algorithm,algorithmic}
\usepackage{lmodern,bookmark,anyfontsize}
\usepackage{booktabs}
\usepackage{multirow}
\usepackage{array} 
\usepackage{hyperref}
\usepackage[capitalise]{cleveref}
\title[Structured column subset selection]{Structured column subset selection with an application to optimal experimental design}

\renewcommand{\t}{^\top}
\newcommand{\B}[1]{\boldsymbol{#1}}
\newcommand{\Bh}[1]{\widehat{\boldsymbol{#1}}}
\newcommand{\T}[1]{\boldsymbol{\mathcal{#1}}}
\newcommand{\gaussianNoise}[1]{\mathcal{N}(\mathbf{0}, \Gnoise)}

\newcommand{\bigO}[1]{\mathcal{O}\left(\displaystyle #1\right)}

\newcommand{\Gnoise}{\B{R}}
\newcommand{\Gpost}{\B{\Gamma}_\textnormal{post}}
\newcommand{\Gprior}{\B{\Gamma}_\textnormal{pr}}
\newcommand{\upost}{\B{u}_\textrm{post}}

\newcommand{\BXj}{\B{X}_{(j)}}

\newcommand{\cost}{\mathcal{C}}

\newcommand{\Th}[1]{\widehat{\boldsymbol{\mathcal{#1}}}}
\newcommand{\R}{\mathbb{R}}

\newcommand{\mc}[1]{\mathcal{#1}}
\newcommand{\logdet}{\mathrm{logdet}}

\newcommand{\bmat}[1]{\begin{bmatrix} #1 \end{bmatrix} }

\newtheorem{lemma}{Lemma}

\usepackage{color}
\usepackage[usenames,dvipsnames]{xcolor}

\title{Structured Column Subset Selection for Bayesian Optimal Experimental Design}
\author[1]{\fnm{Hugo} \sur{D\'iaz}}\email{hsdiazno@ncsu.edu}
\author[1]{\fnm{Arvind } \sur{K.\ Saibaba}}\email{asaibab@ncsu.edu}
\author[2]{\fnm{Srinivas} \sur{Eswar}}\email{seswar@anl.gov}
\author[2]{\fnm{Vishwas} \sur{Rao}}\email{vhebbur@anl.gov}
\author[2]{\fnm{Zichao Wendy} \sur{Di}}\email{wendydi@anl.gov}
\affil[1]{\orgdiv{Department of Mathematics}, \orgname{North Carolina State University}, \orgaddress{\street{2311 Sullivan Drive}, \city{Raleigh}, \state{North Carolina}, \country{USA}, \postcode{27695}}}
\affil[2]{\orgdiv{Mathematics and Computer Science Division}, \orgname{Argonne National Laboratory}, \orgaddress{\street{9700 S. Cass Avenue}, \city{Lemont}, \state{Illinois}, \country{USA}, \postcode{60439}}}
\begin{document}

\abstract{
   We consider optimal experimental design (OED) for Bayesian inverse problems, where the experimental design variables have a certain multiway structure.
   Given $d$ different experimental variables with $m_i$ choices per design variable $1 \le i\le d$, the goal is to select $k_i \le m_i$ experiments per design variable. 
   Previous work has related OED to the column subset selection problem by mapping the design variables to the columns of a matrix $\B{A}$. 
   However, this approach is applicable only to the case $d=1$ in which the columns can be selected independently. 
   We develop an extension to the case where the design variables have a multi-way structure.
   Our approach is to map the matrix $\B{A}$ to a tensor and perform column subset selection on mode unfoldings of the tensor. 
   We develop an algorithmic framework with three different algorithmic templates, and randomized variants of these algorithms. We analyze the computational cost of all the proposed algorithms and also develop greedy versions to facilitate comparisons.  Numerical experiments on four different applications---time-dependent inverse problems, seismic tomography, X-ray tomography, and flow reconstruction---demonstrate the effectiveness and scalability of our methods for structured experimental design in Bayesian inverse problems.
}

\maketitle

\section{Introduction and Motivation}\label{sec:Intro}

Optimal experimental design (OED) aims to enhance inference and decision-making by strategically selecting experiments that maximize a given statistical criterion while respecting physical and budgetary constraints. These techniques are essential in various fields, including nuclear physics \cite{ComptonScattering}, medical imaging \cite{Dale1999, fingerprinting}, and geophysical exploration \cite{seidler2016optimal, geophysicalapplications}. For a comprehensive review, see \cite{Huan_2024}. 

We consider OED in the context of Bayesian inverse problems, where unknown parameters are estimated by combining prior knowledge with data. This approach models parameters, observations, and noise as random variables, with the solution given by a posterior distribution, that is, the conditional distribution of parameters given observations \cite{Dashti2017}. 
This probabilistic framework enables uncertainty quantification and guides data collection through OED \cite{alexanderian2021optimal,eswar2023optimal}.
When the relationship between parameters and data is linear and the noise follows a Gaussian distribution, the posterior distribution remains Gaussian, simplifying inference. 
 In this setting, many common OED criteria admit closed-form expressions, making it possible to evaluate and optimize designs efficiently.

More broadly, OED can be formulated using alternative optimality criteria, depending on the inference objectives and the structure of the problem \cite[Section 4]{alexanderian2021optimal}. For example, A-optimality seeks to reduce the average posterior variance across all parameters, treating all directions in parameter space equally. In contrast, C-optimality targets the variance in a specific linear combination of parameters, focusing the design on improving estimates of a particular quantity of interest.
In this article, we focus on OED based on the D-optimality criterion, which is related to the determinant of the posterior covariance matrix. Equivalently, it can also be interpreted as the expected information gain (EIG) from the prior to the posterior distribution of the parameters of interest \cite{alexanderian2021optimal,Chaloner}.

\paragraph{Structured Column Subset Selection:}
In this paper, we will follow the approach in~\cite{eswar2023optimal} and view OED through the lens of column subset selection.
The column subset selection problem (CSSP) was originally formulated to identify a subset of columns that best preserve the spectral properties of a matrix~\cite{martinsson2020randomized}. In the context of OED, we consider the matrix $\B{A} \in \R^{N \times M}$, which represents a transformed version of the forward model that incorporates prior knowledge and measurement noise \cite{eswar2023optimal}. Details on the construction of this matrix are given in Section~\ref{sec:BOED}. The choice of which columns of $\B{A}$ to select determines which design variables--such as sensor locations, experimental conditions, or control parameters--will be used for inference.
CSSP can be adapted to experimental design by interpreting the columns of a specific matrix $\B{A}$ as candidate designs.
In this context, selecting a subset of columns corresponds to determining design variables to be included in the experiment to maximize informativeness.

In many real-world applications, however, selection is subject to structural constraints that prevent choosing columns independently. 
To capture such constraints, our formulation enforces selection at the level of columns. To make this idea more concrete, we consider the case of time-dependent inverse problems, where there are $m_\text{cs}$ candidate sensor locations  and each sensor can collect $m_{t}$ measurements in time. The goal is to choose a subset of the sensors $K$ (out of $m_\text{cs})$ at which to place the sensors. 
 Specifically, the matrix \( \B{A} \in \R^{N \times M} \) is partitioned as  
\[
\B{A} = \bmat{\B{A}_1 & \cdots & \B{A}_{m_\text{cs}}} \in \R^{N \times M},
\]  
where \( M = m_\text{cs}\cdot m_\text{t} \) and each block \( \B{A}_j \in \R^{N \times m_\text{t}}\) for $1 \le j \le m_\text{cs}$ contains \( m_\text{t} \) columns. 
The goal of OED  is to select \( K \le m_\text{cs} \) representative blocks corresponding to the sensors that are placed. Na\"ively applying column subset selection to select $K\cdot m_\text{t}$ columns does not guarantee that exactly $K$ unique sensors will be selected; it may be possible that the algorithm determines a few columns corresponding to all the sensors.     In other examples, in medical imaging, sensor arrays must be selected as a whole to maintain spatial coherence; and in geophysical monitoring, sensors are deployed in fixed grid patterns, restricting individual placement. Addressing these constraints requires structured selection methods that enforce group-wise selection while retaining the benefits of CSSP. 
More generally, if selection involves \( d \) categorical variables, each with \( m_j \) possible choices in each category, the problem reduces to selecting \( k_j \) elements per category for $1 \le j \le d$. We call this problem a structured column subset selection problem; if \( d=1 \), this formulation recovers the standard CSSP. This is the central focus of the paper.

\paragraph{Challenges and Our Approach:}  
It is well known that CSSP, which is a special case of structured column subset selection, is NP-hard (see, e.g.,~\cite{civril2009selecting}). Enforcing structured selections further restricts the feasible set of solutions. 
Exhaustive search is impractical, and standard greedy heuristics currently do not incorporate structured dependencies.

To address the challenges imposed by the constraints, we propose a structured generalization of CSSP that enforces selection at the level of column blocks, ensuring that structural dependencies are maintained. Our approach leverages tensor decompositions and randomized methods to efficiently capture correlations among design variables while reducing computational costs. By integrating structured selection into CSSP, our methods provide a scalable and effective framework for experimental design in complex inverse problems.

\paragraph{Contributions:}  We summarize the main contributions of this paper:
\begin{enumerate}
\item We propose three novel tensor-based structured selection algorithmic templates---{\it IndSelect}, {\it SeqSelect}, and {\it IterSelect}---that leverage low-rank Tucker-like decompositions and CSSP strategies~\cite{eswar2023optimal,GuEisenstat} to efficiently identify informative subsets of design variables OED for inverse problems (Section~\ref{sec:SSSPAlg}). 

\item We develop a randomized approach (Section~\ref{ssec:rand}) that reduces computational complexity by sketching the matrix $\B{A}$ from the left, to reduce its row dimensions but preserve the column size. We can then use one of the three templates described earlier. The randomized approach also has the benefit that it does not require the adjoint of the forward operator, making it more valuable in many applications. Additionally, it does not require forming the matrix $\B{A}$ explicitly.  
\item We also provide a rigorous computational complexity analysis of all the algorithms (in Section~\ref{ssec:costs}), demonstrating their scalability for large-scale inverse problems. 
\item  We propose greedy versions of the proposed algorithms (in Appendix~\ref{sec:GreedySSSP}) to facilitate comparisons with the proposed algorithms. These extensions also highlight the versatility of the algorithmic templates, since they are readily extensible. 
\end{enumerate}
We validate our methods across four challenging inverse problems—time-dependent partial differential equations (PDEs), seismic imaging, X-ray tomography, and flow reconstruction—demonstrating the broad applicability and effectiveness of our methods in capturing complex dependencies. In numerical experiments, we observe speed ups of up to $50\times$ compared to the corresponding greedy approaches.

\paragraph{Related Work:} 
While a lot of existing work has focused on CSSP (i.e., for $d=1$), we could not find many works that tackle the structured column subset selection. Eswar et al.~\cite{eswar2023optimal} developed approximation algorithms for OED using CSSP algorithms. Our work builds on this idea by extending it to structured data that can be formulated as a tensor, allowing for structured selection across multiple modes of a tensor. A specific generalization for time-dependent problems was considered in~\cite{alexanderian2021optimal}. To our knowledge, none of the approaches in the literature directly and systematically tackle the issue of structured column subset selection, which is the focus of this paper. For this reason, we had to develop greedy versions of the proposed algorithms, to facilitate comparisons between algorithms.   

However,  some works are tangentially related to our goals. In tensor decompositions, subset selection has been used to obtain interpolative decompositions, as in~\cite{saibaba2016hoid,Minster}. Tensor-based sensor placement was also considered in recent work~\cite{Tensor-based_flow} but is different from the proposed approaches in the objective function used. In optimization, subset selection has been enforced across groups of variables using the notion of group sparsity (see, e.g.,~\cite{huang2010benefit}). The notion of clustering different modes of a tensor is the main idea behind a technique called tensor co-clustering~\cite{spectral_co-clustering}.  This framework also shares conceptual similarities with sparse reconstruction techniques in MRI, where prior information is leveraged to recover images from undersampled data \cite{Lustig,Yang2016SparseMRI}.

\section{Preliminaries}
This section provides the necessary notation and background for the discussion that follows. We cover key concepts in tensor theory, Bayesian inverse problems, OED, and column subset selection algorithms,  which form the foundation of the methods and results presented in this article.

\subsection{Tensor Background}
We provide a brief introduction to tensors, focusing on the concepts relevant to this paper; for a more comprehensive treatment, see~\cite{kolda2009tensor}.
Let $\T{X} \in \R^{n_1\times \dots \times n_P}$ be an order-\( P \) tensor with entries \( x_{i_1,\dots,i_P} \), where \( 1 \leq i_\ell \leq n_\ell \) and \( 1 \leq \ell \leq P \).  
Tensors can be unfolded into matrices in different ways, a useful feature for performing linear algebraic operations.
\paragraph{Mode Unfolding}
 A mode-$j$ unfolding (matricization), denoted $\BXj$, reshapes the tensor to form a matrix of size $ {n_j \times  \prod_{\ell=1,\ell \neq j}^P n_\ell}$.
Tensors can be multiplied with matrices using mode products, which are defined via mode unfoldings.
For example, given a ${t \times n_j}$ matrix $\B{S}$, the mode product along mode $j$, denoted, $$\T{Y} = \T{X}  \times_j \B{S} \in \R^{n_1 \times \dots \times n_{j-1} \times t \times n_{j+1} \times \cdots \times n_P},$$  can be computed as  $\B{Y}_{(j)} = \B{S}\B{X}_{(j)}$. 
If  $\B{S}$ and $\B{B}$ are matrices of compatible sizes, mode products satisfy the associativity property  $\T{X}\times_i \B{S}\times_j\B{B} = \T{X} \times_j\B{B} \times_i \B{S}$. For the same mode,  $i= j$,   $\T{X}\times_i \B{S}\times_i\B{B} = \T{X} \times_i\B{BS}$. 

Given two matrices $\B{S}\in \R^{m \times n},\B{B}\in \R^{p \times q}$ ($m,n,p,q$ all positive integers, but with no constraints on the size), then the Kronecker product between the matrices, denoted $\B{S}\otimes \B{B}$, is defined as 
\[ \B{S} \otimes \B{B} \equiv 
\bmat{s_{11} \B{B} & \cdots & s_{1n}\B{B} \\ \vdots & \ddots & \vdots \\ a_{m1}\B{B}  & \cdots & s_{mn} \B{B}} \in \R^{(mp)\times (nq)}. \]
Note that 
$
  \left(\B{S} \otimes \B{B} \right)\t= \B{S}\t \otimes \B{B}\t  \in \R^{ (nq) \times (mp)}.
$
\paragraph{Relation between Mode and Kronecker Products}
Kronecker products are closely related to mode products as follows: Let $\B{S}_j \in \R^{a_j \times n_j}$ for $1\leq j \leq P$ be a sequence of matrices such that the tensor $\T{Y} = \T{X}\bigtimes_{j=1}^P\B{S}_j$ is well-defined. Then 
\begin{align}\label{eq:unfolding_TensorContraction}
    \B{Y}_{(j)} = \B{S}_j \B{X}_{(j)} \left(\B{S}_P \otimes \dots \otimes \B{S}_{j+1} \otimes \B{S}_{j-1} \otimes \dots \otimes \B{S}_1 \right)\t \qquad 1\leq j \leq P.
\end{align}

\paragraph{Tucker Representation and HOSVD}  
Given a tensor \(\T{X} \in \R^{n_1 \times \cdots \times n_P}\) and a target rank \(\B{r} = (r_1, \ldots, r_P)\), the Tucker decomposition approximates \(\T{X}\) as  
\begin{align}\label{def:TuckerFactors}
\T{X} \approx \T{G} \bigtimes_{\ell=1}^P \B{U}_\ell,
\end{align} 
where \(\T{G} \in \R^{r_1 \times \cdots \times r_P}\) is the core tensor and \(\B{U}_\ell \in \R^{n_\ell \times r_\ell}\) are factor matrices with orthonormal columns. 

The higher-order singular value decomposition (HOSVD)~\cite{de_lathauwer_multilinear_2000} generalizes the matrix SVD to tensors by computing the Tucker factor matrices from the truncated SVD: $\B{X}_{(j)}  \approx \B{U}_j \B\Sigma_j\B{V}_j\t$, 
where \(\B{U}_j\)  is a matrix whose columns consist of the leading \(r_j\) singular vectors. The core tensor is then computed as  
\begin{align}\label{def:HOSVDCore}
    \T{G} := \T{X} \bigtimes_{\ell=1}^P \B{U}_\ell\t.
\end{align}  
\paragraph{Sequentially Truncated HOSVD (ST-HOSVD)}
 ST-HOSVD~\cite{vannieuwenhoven_new_2012, Minster} refines HOSVD by computing truncated factor matrices sequentially, mode by mode. 
Rather than optimizing each factor matrix $\B{U}_\ell$ from the full tensor $\T{X}$, it prioritizes preserving the interaction between the core tensor and the factor matrices, as captured in \eqref{def:HOSVDCore}. 
ST-HOSVD
updates the core tensor after each step, ensuring that subsequent truncations are applied to increasingly compressed representations. 
The algorithm proceeds as follows: Initialize the core tensor as $\T{G}^{(0)} = \T{X}$. For each mode $\ell\in\{1,\ldots,P\}$, compute the truncated SVD of the mode-$\ell$ unfolding of \(\T{G}^{(\ell-1)}\),  
\[
\B{G}^{(\ell-1)}_{(\ell)} \approx \widetilde{\B{U}}_{r_\ell} \B{\Sigma}_{r_\ell} \B{V}_{r_\ell}^\top,
\]
retain the top \(r_\ell\) left singular vectors \(\B{U}_\ell=\widetilde{\B{U}}_{r_\ell}\), and update the core tensor as  
\[
\T{G}^{(\ell)} := \T{G}^{(\ell-1)} \bigtimes_{j=1}^{\ell} \B{U}_j^\top \qquad 1 \le \ell \le P.
\]

\subsection{Bayesian Inverse Problems and Optimal Experimental Design} \label{sec:BOED} 
Bayesian inverse problems aim to estimate unknown parameters by integrating prior knowledge with noisy observations within a statistical framework. The solution yields a probability distribution on the unknown parameters, which facilitates uncertainty quantification \cite{alexanderian2021optimal, Dashti2017,Huan_2024}.  
More specifically, a linear inverse problem is formulated as recovering an unknown parameter  $\B{u}\in \R^n$ from observations  $\B{d}\in \R^m$, modeled as
 \begin{align}\label{eq:main}
     \B{d} = \B{Fu}+\B{\varepsilon},
 \end{align}
 where $\B{F}\in \R^{m\times n}$ is the parameter-to-observable map and  $\B{\varepsilon} \sim \gaussianNoise{noise}$ accounts for measurement and model errors. A key challenge in inverse problems is that they are ill-posed, since there may be  insufficient observations (i.e., $m\le n$).
 In order to address these issues, a Bayesian formulation incorporates prior knowledge of the parameter, modeled as  $\B{u} \sim \mathcal{N}(\B{u}_{pr}, \Gprior)$,  leading to the likelihood:
 $$
 \pi_\textrm{like}(\B{d}|\B{u}) \propto \exp\left( -\frac{1}{2}\| \B{Fu}-\B{d}\|^2_{\Gnoise^{-1}}\right).
 $$ 
 With Bayes' rule, the posterior distribution's density $ \pi_\textrm{post}(\B{u}|\B{d})$ satisfies
 \begin{align}\label{def:posterior}
     \pi_\textrm{post}(\B{u}|\B{d}) = 
     \frac{\pi_\textrm{like}(\B{d}|\B{u})\pi_\textrm{pr}(\B{u})} {\pi(\B{d})}\>\propto\> \exp\left( -\frac12\| \B{Fu}-\B{d}\|^2_{\Gnoise^{-1}}
     - \frac12\| \B{u}-\B{u}_\textrm{pr}\|^2_{\Gprior^{-1}}
     \right).
 \end{align}
 Since this problem is linear and the special choice of prior and noise model, the posterior distribution is also Gaussian $\mathcal{N}(\upost,\Gpost)$, where 
 \begin{align}\label{def:post}
   \Gpost^{-1}= \Gprior^{-1}+ \B{F}\t \Gnoise^{-1}\B{F} , ~\mbox{ and }~
     \upost= \B{\Gamma}_\textrm{post}\left(\B{F}\t \Gnoise^{-1}\B{d}
   + \Gprior^{-1}\B{u}_\textrm{pr}
   \right);
 \end{align}
see \cite[Section 3.3]{alexanderian2021optimal}. 
We now cast optimal sensor placement as an OED problem, where a statistical criterion measures optimality.   
\paragraph{Optimality Criterion:}  
We adopt the expected information gain (EIG) as the optimality criterion for sensor selection \cite{alexanderian2021optimal,eswar2023optimal}, which, in finite-dimensional settings, is equivalent to the D-optimality criterion. 
This criterion quantifies the reduction in uncertainty by measuring the volume of the posterior uncertainty ellipsoid, expressed through the determinant of the posterior covariance matrix~\cite[Section 4.1]{alexanderian2021optimal}. It is widely used in Bayesian OED to assess the informativeness of selected observations \cite{Bakker2020ExperimentalDF,eswar2023optimal,TractablGlobal,GOOED}.  
The EIG measures the expected reduction in uncertainty from the prior \(\pi_{\text{pr}}\) to the posterior \(\pi_{\text{post}}\) after observing data. In the Gaussian case, it has a simple analytic form \cite[Section 4.1]{alexanderian2021optimal}, and, up to a multiplicative factor \(1/2\), it is given by  
\begin{align*}
     \phi_{\text{EIG}}:= -\logdet(\Gpost\Gprior^{-1})  
           = \logdet \left(\B{I} + \left(\B{R}^{-\frac 12} \B{F} \Gprior^{\frac 12}\right)\t \left(\B{R}^{-\frac 12} \B{F} \Gprior^{\frac 12}\right)\right).
\end{align*}
Following \cite[Section 3]{eswar2023optimal}, we assume a diagonal noise covariance matrix \(\B{R} = \sigma^2_{\B{R}}\B{I}\), meaning each column of $\Gprior^{\frac 12} \B{F}\t \B{R}^{-\frac 12}$ corresponds to an individual design variable. In contrast, a non-diagonal $\B{R}$ results in linear combinations of the design variables. This suggests defining the matrix \(\B{A}\) as
\begin{align}\label{def:Amatrix}
    \B{A} \equiv \Gprior^{\frac 12} \B{F}\t \B{R}^{-\frac 12}
    = \sigma^{-1}_{\B{R}}\Gprior^{\frac 12} \B{F}\t.
\end{align}
Substituting this into the EIG expression, we get
\begin{align}\label{def:Dopt}
   \phi_\text{EIG}= \logdet(\B{I} + \B{A}\B{A}\t) = \logdet(\B{I} + \B{A}\t\B{A}) = \logdet(\B{I} + \B{\Sigma}_{\B{A}}^2).
\end{align} 
The second equality follows from Sylvester's determinant lemma. Here, $\B{I}$ denotes the identity matrix of appropriate size, and $\B{\Sigma}_{\B{A}}^2$ is a diagonal matrix whose entries are the squared singular values of $\B{A}$, with at most $\operatorname{rank}(\B{A})$ nonzero singular values. We also use the notation $\Psi(\B{A}) = \logdet(\B{I}+\B{AA}\t)$, so that $\phi_\text{EIG} = \Psi(\B{A})$. In the sequel, we show how the design enters the design criterion.

\subsection{Selection Matrices and Subsampled EIG} \label{sec:SelectionMatrices}  

To formalize the selection of a subset of design variables, we introduce a selection matrix \(\B{S}\), which extracts a subset of columns from \(\B{A} \in \mathbb{R}^{N \times M}\). The discussion below corresponds to the case $d=1$, i.e., there is a single design variable.

Let \(\mathcal{S} = \{ {i_1},\ldots, {i_K}\}\) be an index set representing the selected \(K\le M\) columns of \( \B{A}  \). Then the corresponding selection matrix $\B{S}$ is defined as
\begin{align}
    \B{S} = [\B{e}_{{i_1}}, \B{e}_{{i_2}}, \dots, \B{e}_{{i_k}}] 
\in \mathbb{R}^{M \times K},
\end{align}
where \(\B{e}_{i}\) denotes the \(i\)th standard basis vector in \(\mathbb{R}^M\). 
The selected columns of $\B{A}$ are given by
\begin{align}
    \B{A}(:,\mathcal{S})= \B{AS} \in \mathbb{R}^{N \times K}.
\end{align}
We now express the posterior distribution and the corresponding EIG based on a given selection. Note that the EIG is independent of the specific ordering of selected columns.

By restricting the  data to that collected by the selected subset of design variables, we obtain the following expressions for the posterior mean and covariance:
\begin{align}\label{def:postRestricted}
   \B{\Gamma}_\textrm{post}(\mc{S}) &= (\B{\Gamma}_\textrm{pr}^{-1} + \sigma_{\B{R}}^{-2} \B{F}^\top \B{S} \B{S}^\top \B{F})^{-1},  \\
   \B{u}_\textrm{post}(\mc{S}) &= \B{\Gamma}_\textrm{post}({S}) \left( \sigma_{\B{R}}^{-2} \B{F}\t \B{S} \B{S}\t \B{d} + \B{\Gamma}_\textrm{pr}^{-1} \B{u}_\textrm{pr} \right).
\end{align}  

The EIG associated with the selected subset of columns from $\B{A}$ is given by 
$
     \logdet \left( \B{\Gamma}_\textrm{pr}^{\frac 12} \B{\Gamma}_\textrm{post}^{-1}(\mc{S}) \B{\Gamma}_\textrm{pr}^{\frac 12} \right).$
Then,  substituting \eqref{def:postRestricted}, we obtain
\begin{align}\label{def:SubsampledEIG}
    \phi_\text{EIG}(\mc{S} ) :=
     \logdet(\B{I} + \B{AS}(\B{AS})\t).
\end{align}

\subsection{Column Subset Selection Problem}  
\label{sec:CSSPGKS}  

The column subset selection problem aims to identify a subset of columns from a given matrix that preserves key spectral properties \cite{mahoney2009cur,martinsson2020randomized}. This problem is fundamental in dimensionality reduction and low-rank approximations.   
The discussion below corresponds to the case $d=1$. 

\paragraph{Golub--Klema--Stewart Method (GKS)}  
In this work we focus on the GKS method \cite[Section 5.5.7]{MatrixComputations}, a two-stage approach for selecting \( K \) informative columns (\( K \le M \)) from a matrix $\B{A} \in \mathbb{R}^{N \times M}$. The GKS method provides a way to extract a well-conditioned subset of columns that retains essential spectral characteristics of \( \B{A} \).  
First, a rank-\( K \) approximation of \( \B{A} \) is obtained via truncated SVD:  
\[
\B{A} \approx \B{U}_K \B{\Sigma}_K \B{V}_K\t,  
\]
where 
\( \B{U}_K \in \mathbb{R}^{N \times K} \), 
\( \B{\Sigma}_K \in \mathbb{R}^{K \times K} \),
and \( \B{V}_K \in \mathbb{R}^{M \times K} \).  
In the second stage, a column-pivoted QR (CPQR) factorization  of  \(\B{V}_K\t \in \R^{K \times M}\) produces 
\begin{align}\label{eq:CPQRVH}
    \B{V}_K\t \B{\Pi}= \B{V}_K\t \bmat{ \B{\Pi}_1 &  \B{\Pi}_2} = 
    \B{Q} \bmat{\B{R}_{11} &  \B{R}_{12}} 
    =\bmat{\B{V}_{11}\t & \B{V}_{12}\t},
\end{align}
where \( \B{\Pi} \) is a permutation matrix. The first \( K \) columns of \( \B{\Pi} \), denoted \( \B{\Pi}_1 \), determine the selected indices. The corresponding selection matrix is given by \( \B{S} = \B{\Pi}_1 \in \mathbb{R}^{M \times K} \), and the final subset of selected columns is \( \B{W} = \B{A} \B{S} \).  
\Cref{alg:detcssp} provides a formal outline of this method.

In the second stage, any column selection algorithm can be applied, such as strong rank-revealing QR (sRRQR) \cite{GuEisenstat}, DEIM point selection algorithm \cite[Algorithm 1]{DEIMCUR},  leverage scores \cite{MAL-035}, greedy algorithms \cite{Altschuler}, or randomized QR with column pivoting \cite{Duersch, epperly2024embrace}. We choose CPQR because of its strong empirical performance and availability in standard numerical linear algebra libraries.
\paragraph{Spectral Properties of the GKS Method}
For any column selection algorithm used in the second stage of GKS, as long as $\B{V}\t_K\B{S}$ is not singular (see \eqref{eq:CPQRVH}),
 the singular values of the selected columns $\B{AS}$ satisfy the following bound:
\begin{align}\label{def:GKSSpectralBound}
\frac{\sigma_j(\B{A})}{\|(\B{S}\t\B{V}_K)^{-1} \|_2} \le \sigma_j(\B{AS}) \le \sigma_j(\B{A}) \qquad 1 \le j \le K. 
\end{align}

A key consequence of \eqref{def:GKSSpectralBound} is the following lemma, which establishes bounds on the EIG for the selection obtained via GKS. 
\begin{lemma}\label{lemma:csspoed}
Let \( \B{A} \in \mathbb{R}^{N \times M} \) with \( K \leq \operatorname{rank}(\B{A}) \), and let \( \B{S} \) be a selection matrix, corresponding to an index set $\mc{S}$  of \( K \) distinct columns, such that \( \B{V}_K^\top \B{S} \) is nonsingular. Then  
\begin{align}
    \Psi\left(\frac{\B{\Sigma}_K}{\|(\B{S}^\top \B{V}_K)^{-1}\|_2}\right) \leq \phi_\text{EIG}( \mc{S}) \leq \phi_\text{EIG}( \mc{S}^\text{opt}) \leq \Psi(\B{\Sigma}_K) \leq \Psi(\B{A}),
\end{align}
where \( \mc{S}^\text{opt} \) denotes a selection of indices that chooses an optimal set of \( k \)  columns maximizing the EIG.   
\end{lemma}

\begin{proof}
See \cite[Theorem 3.2]{eswar2023optimal}. 
\end{proof}

This result provides theoretical guarantees on the quality of the selected columns, highlighting the effectiveness of the GKS method in preserving spectral properties relevant to the EIG.  
A key advantage of some column selection algorithms is that they provide bounds of the form $1 \le \|(\B{S}^\top \B{V}_K)^{-1}\|_2 \le q_K $. 
 Here \(q_K\) is an upper bound that depends on the particular choice of the algorithm. The smaller $q_K$ can be made, the closer the selected columns are to the optimal.   Table~\ref{tab:bounds} lists several methods, the corresponding value of $q_K$, and the computational complexity of the selection.
\begin{table}[!ht]
    \centering
    \renewcommand{\arraystretch}{1.3} 
    \begin{tabular}{|c|c|c|c|}
        \hline
        \textbf{Method} & \textbf{} \( q_K \) &  \textbf{Comput. Complex.} & \textbf{Reference} \\ 
        \hline
        CPQR/ Q-DEIM  & \( 2^K\sqrt{M-K} \) & \( \mathcal{O}(MK^2+K^3) \) & \cite{GuEisenstat} \\ 
        sRRQR & \( \sqrt{1+f^2K(M-K)} \) & \( \mathcal{O}(K^2 M \log_f(M)) \) & \cite{GuEisenstat} \\ 
        DEIM  & \( 2^K\sqrt{\frac{MK}{3}}\) & \( \mathcal{O}(MK) \) & \cite[Section 2]{DEIMCUR} \\ 
        \hline
    \end{tabular}
    \caption{Bounds and computational complexities for different selection methods. Here $f> 1$ is a user-defined parameter in the sRRQR algorithm. }
    \label{tab:bounds}
\end{table}

\begin{algorithm}[!ht]
\caption{\textsc{DetCSSP}: Deterministic Column Subset Selection via GKS}
\label{alg:detcssp}
\begin{algorithmic}[1]
    \REQUIRE Matrix \(\B{A} \in \mathbb{R}^{N \times M}\), target rank \(K \leq \min\{M, N\}\)
    \ENSURE Selection operator \(\B{S} \in \mathbb{R}^{M \times K}\), submatrix  \(\B{W} \in \mathbb{R}^{N \times K}\)
    
    \STATE Compute the top-\(K\) right singular vectors of \(\B{A}\), denoted by \(\B{V}_K \in \mathbb{R}^{M \times K}\)
    \STATE Perform pivoted QR factorization of \(\B{V}_K^\top\):
    \[
    \B{V}_K^\top 
    \begin{bmatrix} \B{\Pi}_1 & \B{\Pi}_2 \end{bmatrix}
    = \B{Q}
    \begin{bmatrix}
    \B{R}_{11} & \B{R}_{12}
    \end{bmatrix}
    \]
    where \(\B{\Pi}_1 \in \mathbb{R}^{M \times K}\) is the permutation selecting the top \(K\) columns.
    \STATE Set \(\B{S} = \B{\Pi}_1\) and \(\B{W} = \B{A} \B{S} \in \mathbb{R}^{N \times K}\)
    
    \RETURN \(\B{S}, \B{W}\)
\end{algorithmic}
\end{algorithm}

We summarize the CSSP algorithm based on the GKS method in Algorithm~\ref{alg:detcssp}, which consists of computing a truncated SVD followed by a column-pivoted QR factorization on the top right singular vectors.

\paragraph{Computational Complexity of the GKS Method}

This method  requires a truncated SVD on $\B{A}\in \R^{N\times M}$, which has a computational complexity of \(\mathcal{O}( \min(N,M)^2 \cdot \max(N,M))\) flops if a thin SVD is performed.
Additionally, it involves a column-pivoted QR factorization on a short and wide matrix \(\B{V}\t_K \in \R^{K \times M}\), with a complexity of \(2MK^2-2K^3/3\) flops \cite[Section 3.3]{eswar2023optimal}. Thus, the dominant computational cost arises from the truncated SVD, motivating the need for more suitable alternatives for large-scale problems. To address this computational cost, we use randomized approaches in Section~\ref{ssec:rand}.

\section{Algorithms for Structured OED }

This section introduces efficient and scalable approaches for OED with structured column selection. The proposed methods combine tensor decompositions, CSSP techniques,  randomized sketching, and adjoint-free variants.    
Since mode unfolding arranges data row-wise, the problem reduces to row subset selection. In numerical linear algebra, however, column selection is the primary focus.  
To leverage standard tensor and matrix libraries, we first compute the Tucker-like factors for the target modes and then apply a column selection algorithm to their transposes. In Section~\ref{ssec:overview}, we give the overview of our approach and discuss the three algoritmic templates in Section~\ref{sec:SSSPAlg}. In Section~\ref{ssec:rand}, we give the Sketch-First randomized approach and, finally, in Section~\ref{ssec:costs}, we discuss the computational costs of all the proposed algorithms. The algorithmic templates can be adapted to give greedy variants of the algorithms, which we discuss in Appendix~\ref{sec:GreedySSSP}, including a detailed comparison with the proposed algorithms.

\subsection{Overview of Our Approach}  \label{ssec:overview}

We consider tensors derived from the matrix $\B{A} \in \R^{N\times M}$, which corresponds to a weighted forward operator as defined in \eqref{def:Dopt}. We reshape the matrix $\B{A}$ to obtain the tensor $\T{X} \in \R^{m_1 \times \cdots \times m_d \times N}$, where $d$ represents the number of categories of selection variables. The matrix $\B{A}$ can be recovered from the tensor $\T{X}$ as  
\begin{align}\label{def:A}
    \B{A} = \B{X}_{(d+1)} \in \R^{N \times M}, \quad M = \prod_{j=1}^{d} m_j.
\end{align}

The structure of $\T{X}$ enables structured selection across multiple experimental variables.
For each mode $j \in \{1, \dots, d\}$, we select $k_j$ indices, represented by the tuple $\B{k} = (k_1, \dots, k_d)$, where $k_j \leq m_j$ for $1\le j \le d$. The total number of selected design variables  is  
$K = \prod_{j=1}^{d} k_j.$ The following discussion generalizes the approach for $d=1$ in Section~\ref{sec:SelectionMatrices}.
\paragraph{Structured Selection (several design variables)} Given a positive integer $t$,  denote $[t] = \{1,\dots,t\}$. Given $d$ categorical variables, the columns of $\B{A}$ can be identified with the set $[m_1] \times \dots \times [m_d]$. 
We must choose $k_j$ indices out of $[m_j]$ indices for $1\le j \le d$. We denote the sets $\mc{S}^{(j)} = \{i_1^{(j)},\dots, i_{k_j}^{(j)}\}$ so that the selection operator is of the form $\mc{S} = \mc{S}^{(1)} \times \dots  \times \mc{S}^{(d)}$. The corresponding selection operator is 
\[\B{S} = \B{S}_d \otimes \dots \otimes \B{S}_1, \qquad \text{where}\qquad \B{S}_j = \B{I}(:, \mc{S}^{(j)}), \> 1 \le j \le d.  \] 
The corresponding selected columns of $\B{A}$ are denoted $\B{AS} \in \R^{N\times K}$, where $K = \prod_{j=1}^dk_j$.
We can also compute the subsampled tensor as 
\[ \Th{X} = \T{X} \times_1 \B{S}_1\t \dots \times_d \B{S}_d\t.  \]
With this notation $\B{AS} = \Bh{X}_{(d+1)} = \B{X}_{(d+1)} (\B{S}_d \otimes \dots \otimes \B{S}_1)$, 
where $\B{S} = (\B{S}_d \otimes \dots \otimes \B{S}_1)$. The main objective of this article is to find a selection operator \( \B{S} = \B{S}_d \otimes \dots \otimes \B{S}_1 \) that maximizes~\eqref{def:SubsampledEIG}. Since this problem is NP-hard, prior work \cite{eswar2023optimal} has proposed approximation methods based on CSSP (i.e., $d=1$), which selects a subset of columns that approximately preserve the main spectral properties of \( \B{A} \), to identify the most informative design variables. 
We extend this approach to the more general case $d > 1$.

The following methods adapt strategies for computing the Tucker factors, treating the number of selected variables per mode as the target rank. We first discuss the methods and then discuss the computational cost of each method.

\subsection{Algorithmic Templates for Structured OED}\label{sec:SSSPAlg}
We begin by presenting a framework for selecting indices in each mode, ensuring a structured selection, with three different methods.
Each method outputs $d$ selection matrices \((\B{S}_1,\ldots, \B{S}_d)\). 
The methods differ in how they incorporate information from other design variables (modes).
While the proposed framework is compatible with any CSSP algorithm, we primarily employ the GKS-based approach detailed in \Cref{alg:detcssp}. In order to emphasize the fact that any CSSP algorithm can be used (hence, an algorithmic template), the proposed algorithms refer to a function $\B{S} \leftarrow \text{CSSP}(\B{M},k)$, which takes in a matrix $\B{M} \in \R^{m\times n}$ and returns a selection operator $\B{S}$ corresponding to $k \le n$ important columns $\B{MS}$. 

\subsubsection{IndSelect: Independent Mode Selection}  This method selects indices from each mode independently by applying an appropriate CSSP method to the mode unfolding 
\begin{align}\label{def:unfoldInd}
\B{X}_{(j)}\t \in \R^{\left(NM / m_j\right) \times m_j}, \quad 1 \le j \le d.
\end{align}
This is summarized in Algorithm~\ref{alg:indselect}. This  tensor-based OED method resembles the HOSVD for computing a low-rank Tucker decomposition. 

\begin{algorithm}[!ht]
\caption{\textbf{IndSelect: Independent Mode Selection}}
\label{alg:indselect}
\begin{algorithmic}[1]
\REQUIRE Matrix $\B{A}\in \R^{N\times M}$, number of indices \( \B{k} = (k_1, \dots, k_d) \)
\ENSURE Selected indices \( (\B{S}_1, \dots, \B{S}_d) \) for each mode  
\STATE Reshape matrix $\B{A}$ into tensor \( \T{X}  \in \R^{m_1 \times \cdots \times m_d \times N}\)
\FOR{$j=1,\dots,d$ }  
    \STATE \( \B{S}_{j} \gets \text{CSSP}(\B{X}\t_{(j)}, k_j) \)   
    \hfill \texttt{\small \# e.g., apply GKS method:\Cref{alg:detcssp} } 
\ENDFOR  
\RETURN \( (\B{S}_1, \dots, \B{S}_d) \)  
\end{algorithmic}
\end{algorithm}

\subsubsection{SeqSelect: Sequential Mode Selection} 
In this approach, we start with the same as in IndSelect and obtain the selection operator $\B{S}_1$ from $\B{X}_{(1)}\t$. However, in the second step, rather than working with $\B{X}_{(2)}\t$, we form the core tensor $\T{G} = \T{X}\times_1 \B{S}_1\t \in \R^{k_1 \times m_2 \dots \times m_d \times N} $ and perform mode subset selection on $\B{G}_{(2)}$. We can thus proceed sequentially. 

At core $j$ (for $1\le j \le d$), we can consider the core tensor 
\begin{align}\label{def:unfoldSeq}
\T{G}^{(j-1)} = \T{X} \bigtimes_{i=1}^{j-1}\B{S}_i\t \in \R^{k_1 \dots \times k_{j-1}\times m_j  \dots \times m_d \times N},
\end{align}
with $\T{G}^{(0)} = \T{X}$, and we perform subset selection on $\B{G}^{(j-1)}_{(j)}$ to obtain $\B{S}_j$. Thus, this method incorporates the selection  performed in the previous steps. It is similar to the ST-HOSVD approach for computing low-rank Tucker decompositions.  Similar to ST-HOSVD, we allow for a different order of processing the modes; however, unlike ST-HOSVD, we  focus only on reducing the core tensor and do not consider the factor matrices. This approach is summarized in Algorithm~\ref{alg:seqselect}. 

\begin{algorithm}[!ht]
\caption{\textbf{SeqSelect:  Sequential Mode Selection}}
\label{alg:seqselect}
\begin{algorithmic}[1]
\REQUIRE Matrix $\B{A} \in \R^{N\times M}$, number of indices \( \B{k} = (k_1, \dots, k_d) \), processing order \( \B{\rho} = (\pi_1, \dots, \pi_d) \)  
\ENSURE Selected indices \( (\B{S}_1, \dots, \B{S}_d) \)  
\STATE Reshape matrix $\B{A}$ as tensor \( \T{X}  \in \mathbb{R}^{m_1 \times \cdots \times m_d \times N} \) 
\STATE \( \T{G} \gets \T{X} \) \hfill \texttt{\small \# Create a copy of \(\T{X}\) to preserve the original tensor}  
\FOR{$j=1,\dots,d$}  
    \STATE \( \B{S}_{j} \gets 
    \text{CSSP}(\B{G}_{(\pi_j)}\t, k_{\pi_j}) \)  
    \STATE \( \B{G}_{(\pi_j)} \gets \B{S}_{\pi_j}\t \B{G}_{(\pi_j)} \) \hfill \texttt{\small \# Updates $\T{G}$, keeping only $k_{\pi_j}$ rows from current unfolding}  
\ENDFOR  
\RETURN \( (\B{S}_1, \dots, \B{S}_d) \)  
\end{algorithmic}
\end{algorithm}

\subsubsection{IterSelect: Iterative Mode Selection} 
 IterSelect begins with a random initialization of the selection operators defined by $(\B{S}_1^{(0)},\dots,\B{S}_d^{(0)}).$ Then the method performs a sweep through all the modes to sequentially determine the selection operator. At iteration $t$, we consider mode $j$ ($1 \le j \le d$). We have the selection operator defined by $$\left(\B{S}_1^{(t)}, \dots,\B{S}_{j-1}^{(t)}, \B{S}_{j}^{(t-1)}, \dots,\B{S}_d^{(t-1)}\right).$$
To obtain the selection operator $\B{S}_j^{(t)}$, we form the intermediate tensor 
\begin{align}\label{def:unfoldIter}
\T{Y}^{(t,j)} = \T{X} \bigtimes_{i=1}^{j-1} \left(\B{S}_i^{(t)}\right)\t \bigtimes_{i=j+1}^d \left(\B{S}_i^{(t-1)}\right)\t \in \R^{k_1 \times \dots\times k_{j-1} \times m_j \times k_{j+1}\times \dots \times k_d \times N},
\end{align}
for $t \ge 0$ and  $1 \le j \le d$. At the end of each sweep through the modes, 
the method checks for improvement in the design criterion $\phi_\text{EIG}$. 
If the improvement is below a specified tolerance or if EIG decreases, the algorithm terminates and returns the last selection matrices.
This method is similar to the higher-order orthogonal iteration (HOOI) \cite{BestRank1,de_lathauwer_multilinear_2000, HOOI} 
method for Tucker approximation and is summarized in Algorithm~\ref{alg:iterselect}. It also has similarities to block coordinate descent approaches~\cite{beck2013convergence}.

\begin{algorithm}[!ht]
\caption{\textbf{IterSelect: Iterative Mode Selection}}
\label{alg:iterselect}
\begin{algorithmic}[1]
\REQUIRE Matrix $\B{A}\in \R^{N\times M}$, number of indices \( \B{k} = (k_1, \dots, k_d) \), processing order \( \B{\rho} = (\pi_1, \dots, \pi_d) \),  tolerance \texttt{tol}
\ENSURE Selected indices \( (\B{S}_1, \dots, \B{S}_d) \)  
\STATE Reshape matrix $\B{A}$ into tensor \( \T{X}  \in \mathbb{R}^{m_1 \times \cdots \times m_d \times N} \)
\STATE {Initialize:} Select \( \B{S}_1, \dots, \B{S}_d \) randomly and define $\B{S} = \B{S}_d\otimes \dots \otimes \B{S}_1$  

\STATE {Set initial EIG:} \( \phi_{\text{prev}} \gets -\infty \), $\phi_{\text{curr}} \gets \phi_\text{EIG}(\B{S})$

\WHILE{\textbf{true}}  
    \STATE \( \widetilde{\B{S}}\gets \B{S}
    \)
    \FOR{ \( j = 1,\dots,d \)}  
        \STATE {Reset selection matrix:} \( \widetilde{\B{S}}_{\pi_j} \gets \B{I}_{m_{\pi_j}} \)  
        \STATE {Apply selection:}  
        \(
        \T{Y} \gets \T{X} \bigtimes_{j=1}^{d} \widetilde{\B{S}}_j\t
        \)   
        \STATE {Update selection:}  
        \(
        \widetilde{\B{S}}_{\pi_j} \gets \text{CSSP}\left(\B{Y}_{(\pi_j)}\t, k_{\pi_j}\right)
        \) 
    \ENDFOR  
   \STATE Compute \( \phi_{\text{curr}} \gets \phi_\text{EIG}(\widetilde{\B{S}}) \)
       \IF{\( |\phi_{\text{curr}} - \phi_{\text{prev}}|/\phi_{\text{prev}} < \texttt{tol} \textbf{ or } \phi_{\text{curr}} < \phi_{\text{prev}}
    \)
    }
        \STATE \textbf{break}
    \ENDIF
    \STATE \( {\B{S}}\gets \widetilde{\B{S}}, \phi_{\text{prev}} \gets \phi_{\text{curr}} \)
\ENDWHILE  

\RETURN \( (\B{S}_1, \dots, \B{S}_d) \)  
\end{algorithmic}
\end{algorithm}

After completing a full sweep through all modes, the algorithm evaluates the  EIG  criterion. If the improvement between two consecutive iterations falls below a specified tolerance, here set to $10^{-10}$, or if the criterion value decreases (which may occur due to lack of monotonicity guarantees), the iteration stops.  As in SeqSelect, we can prescribe a processing order for the modes. Furthermore, instead of a random initialization, we can consider the outputs of IndSelect and SeqSelect as initial guesses for IterSelect; we do not consider it in this paper.

\subsection{Randomized Approach: Sketch-First, then Structured Column Selection}\label{ssec:rand}
A key bottleneck in the GKS method is the computation of the truncated SVD. To mitigate this, we introduce a \textit{Sketch-First} approach that leverages randomized sketching to reduce the size of the matrix before applying structured subset selection. 

\begin{figure}[h]
    \centering
    \includegraphics[width=0.8\textwidth]{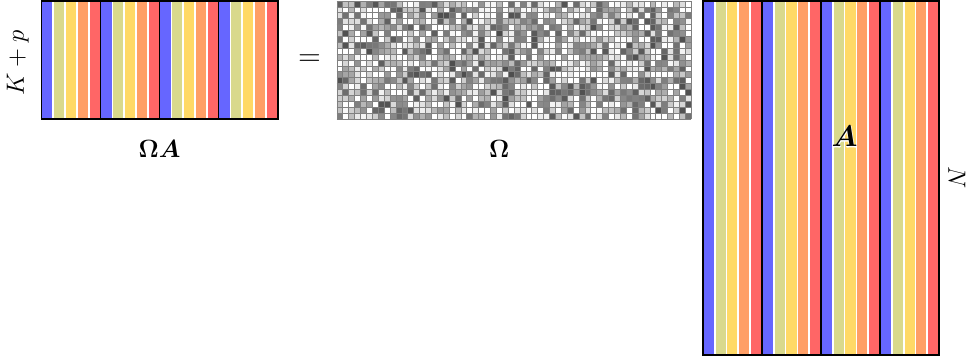}
    \caption{Sketching process: matrix \( \B{A} \in \mathbb{R}^{N \times M} \) is ``projected" via a random matrix \( \B{\Omega} \in \mathbb{R}^{r \times N} \) to obtain a compressed (row-wise) matrix \( \B{Y} = \B{\Omega} \B{A} \in \mathbb{R}^{r \times M} \).}
    \label{fig:sketching}
\end{figure}

In this approach we draw a random matrix  \( \B{\Omega} \in \R^{r \times N} \) with \(r=  p+\prod_{j=1}^d k_j \) rows, where entries are drawn from a Gaussian distribution, \( \B{\Omega} \sim \mathcal{N}(0, 1/r) \). Here \( p \)  is an oversampling parameter typically taken to be $\le 20$. Then, we form the sketch  \( \B{Y} = \B{\Omega} \B{A} \in \R^{r \times M} \), which is much smaller than $\B{A}$ if we assume that $K  \ll N$; see \Cref{fig:sketching}.  
This assumption is met in all the numerical experiments we consider. The sketched matrix $\B{Y}$ (approximately) preserves the largest singular values of $\B{A}$.

The main idea is to apply the selection methods---IndSelect, SeqSelect, and IterSelect---to the sketch matrix $\B{Y}$ rather than $\B{A}$. In particular, we can reshape $\B{Y}$  
into a tensor $\T{Y}$ of size $m_1 \times \dots \times m_d \times r$. In Algorithm~\ref{alg:indselect} we show how to combine the sketching with IndSelect. For brevity, we omit the other two variants involving SeqSelect and IterSelect.

\begin{algorithm}[!ht]
\caption{\textbf{Sketch-First-IndSelect}}
\label{def:algorithmSketchedFirstIndSelect}
\begin{algorithmic}[1]
\REQUIRE Matrix \( \B{A} \in \R^{N \times M} \),  number of indices \( \B{k} = (k_1, \dots, k_d) \),  oversampling parameter \( p \) 
\ENSURE Selected indices \( (\B{S}_1, \dots, \B{S}_d) \) for each mode  

\STATE Compute \(r=p+K \) with \( K = \prod_{j=1}^d k_j\)
\STATE Generate a random matrix \( \B{\Omega} \in \R^{r \times N} \) with independent entries drawn from $\mathcal{N}(0,1/r)$ 
 
\STATE Compute sketched matrix \( \B{Y} = \B{\Omega} \B{A} \)  

\STATE Reshape \( \B{Y} \) into tensor \( \T{Y} \in \R^{m_1 \times \cdots \times m_d \times r} \)  

\FOR{ $j=1,\dots,d$}  
    \STATE \( \B{S}_{j} \gets \text{CSSP}(\B{Y}\t_{(j)}, k_j) \)   
    \hfill \texttt{\small \# Apply GKS or another selection method} 
\ENDFOR 
\RETURN \( (\B{S}_1, \dots, \B{S}_d) \)  
\end{algorithmic}
\end{algorithm}

The method IterSelect requires computing the EIG  \( \logdet(\B{I} + \B{A}\B{S}(\B{A}\B{S})\t) \). In the sketched setting, we instead evaluate \( \logdet(\B{I} + (\B{YS})(\B{Y}\B{S})\t) \).
This substitution is justified by the fact that \(\B{Y} = \B{\Omega} \B{A} \) approximately preserves the dominant spectral properties of \( \B{A} \), particularly its leading singular values. See \cite[Section 4]{eswar2023optimal} for more precise statements for $d=1$.

\paragraph{Benefits} Besides computational efficiency, the proposed method has other benefits. First, observe that IndSelect, SeqSelect, and IterSelect all require the matrix $\B{A}$ to be formed since it needs to be reshaped into a tensor. By first sketching, however, we eliminate the need to form $\B{A}$ explicitly. To form the sketch, we can compute 
\begin{align}\label{def:OmegaA}
    \B{Y}\t = \sigma_{\B{R}}^{-1} \B{F} (\B\Gamma_{\rm pr}^{1/2}\B\Omega\t) \in \R^{M \times r}.
\end{align}
We note that the action of \( \B{F}^\top \) is eliminated and only the action of $\B{F}$ is needed. This is especially beneficial since the adjoint operator $\B{F}\t$ is unavailable (due to legacy software) or expensive to apply in many applications. Furthermore, the sketch can be parallelized by applying $\B\Gamma_{\rm pr}^{1/2}$ followed by the forward operator, in parallel, across each column of $\B\Omega\t$. Thus, this method can be viewed as an extension of the RAF-OED approach from \cite[Section 4]{eswar2023optimal} to the structured selection case.

Other ways exist to incorporate randomization in the structured column selection. We also test a randomized GKS variant by replacing the SVD with randomized SVD (for example, ~\cite[Algorithm 4.4]{HalkoMartinssonTropp}) in each unfolding \( \B{X}_{(j)}^\top \). However, sketching each unfolding separately increases computational cost. At the same time, unlike the Sketch-First approach, it also requires forming $\B{A}$ explicitly. Therefore, we do not pursue these approaches further.

\subsection{Computational Cost}  \label{ssec:costs}
We analyze the computational cost of the proposed methods. We emphasize that the algorithms \textit{IndSelect}, \textit{SeqSelect}, and \textit{IterSelect} can be applied with any CSSP method, not just GKS. Therefore, in the cost analysis, we first present the cost for a generic CSSP method and then discuss the cost with GKS. To this end, 
let \( \cost_\text{cssp}(n_\text{c}, n_\text{r}, k) \) denote the computational cost  of selecting \(k\) columns from a matrix of size \(n_\text{r} \times n_\text{c}\). For example, the cost of GKS (a truncated SVD followed by pivoted QR) is $\cost_\text{GKS}(n_\text{c}, n_\text{r}, k) = \mc{O}(n_cn_r\min\{n_c,n_r\} + n_ck^2)$ flops.  If RandGKS is used (randomized SVD followed by pivoted QR), then the cost is $\mc{O}(n_cn_rk + n_ck^2)$ flops.

\paragraph{Cost of Forming \texorpdfstring{$\B{A}$}{} and Sketching} The input to the main algorithms requires $\B{A}$ as an input. Therefore, in discussing the computational cost of the approaches, we must factor in the cost of forming $\B{A}$. In the applications we consider, we can access $\B{A}$ through matrix-vector products, or matvecs, involving $\B{A}$ and its transpose. We assume that  \( T_{\B{A}} \) denotes the cost of applying \( \B{A} \) to a vector and \( T_{\B{A}^\top} \) is the cost of applying \( \B{A}^\top \). 
The matrix $\B{A}$ can be formed by multiplying by the identity matrix from the left or right. On the other hand, the cost of forming $\B{Y} = \B\Omega\B{A}$ requires $r$ matvecs with $\B{A}\t$. 
Table~\ref{tab:cost_sketch_vs_deterministic} summarizes the cost of constructing these matrices.

\begin{table}[!ht]
\centering
\renewcommand{\arraystretch}{1.4}
\begin{tabular}{@{}lcc@{}}
\toprule
\textbf{Approach} & \textbf{Matrix Constructed} & \textbf{Computational Cost} \\
\midrule
Deterministic & \( \B{A} \in \mathbb{R}^{N \times M} \) & \( \min\left\{ M \cdot T_{\B{A}},\; N \cdot T_{\B{A}^\top} \right\} \) \\
Sketch-First        & \( \B{\Omega A} \in \mathbb{R}^{(K+p) \times M} \) & \( (K + p) \cdot T_{\B{A}^\top} \) \\
\bottomrule
\end{tabular}
\captionsetup{width=\textwidth}
\caption{Comparison of the computational cost of forming  the matrix \( \B{A} \) (deterministic approach) versus  \( \B{\Omega A} \) (Sketch-First approach). Here, \( T_{\B{A}} \) denotes the cost of applying \( \B{A} \) to a vector, and \( T_{\B{A}^\top} \) is the cost of applying \( \B{A}^\top \). }
\label{tab:cost_sketch_vs_deterministic} 
\end{table}

Once the input matrix is formed, the remaining cost comes from applying a CSSP method to each selected mode unfolding. 
This cost depends on the size of each  unfolding and the CSSP method used.
We assume the input matrix is \( \B{A} \in \mathbb{R}^{N \times M} \), which leads to a tensor $\T{X}  \in \mathbb{R}^{m_1 \times \cdots \times m_d \times N}$; the cost for the Sketch-First approach is discussed later. 
In what follows, \( k_j \) is the number of selected rows in the \( j \)-th mode (i.e., design variables), and the total number of design configurations is \( K = k_1 \cdots k_d \).

\paragraph{IndSelect}  
This method applies a CSSP method independently to each mode unfolding of \( \T{X} \). The transpose of the \( \ell \)th unfolding has shape \( (NM / m_\ell) \times m_\ell \); see \eqref{def:unfoldInd}. Thus, the total cost of this method is
\[
\sum_{\ell=1}^{d} \cost_\text{cssp} \left(m_\ell, N(M/{m_\ell}), k_\ell \right) \quad \text{flops}.
\]

\paragraph{SeqSelect} 
This method reduces the tensor size at each step based on previous selections. As shown in \eqref{def:unfoldSeq}, the transpose of the \( \ell \)th unfolding of \( \T{G}^{(\ell-1)} \) transposed has shape
\[
(k_1 \cdots  k_{\ell-1}) \cdot ( m_{\ell+1}\cdots m_d) \cdot  N \times m_\ell = N\frac{M}{m_\ell} \prod_{i=1}^{\ell-1} \frac{k_i}{m_i} \times m_\ell.
\]
We apply a CSSP approach on the transpose of $\T{G}^{(\ell-1)}, $ which leads to the total cost
\[
\sum_{\ell=1}^{d} \cost_\text{cssp} \left(m_\ell, N\frac{M}{m_\ell} \prod_{i=1}^{\ell-1} \frac{k_i}{m_i}, k_\ell \right) \> \text{flops}.
\]

\begin{table}[t]
\centering
\renewcommand{\arraystretch}{1.6}
\begin{tabular}{@{}lll@{}}
& \textbf{General CSSP Cost} & \textbf{GKS-Based Cost} \\
\midrule
\textbf{IndSelect} &
\( \displaystyle \sum_{\ell=1}^{d} \cost_\text{cssp} \left(m_\ell, N_\text{rows} \cdot \frac{M}{m_\ell}, k_\ell \right) \) &
\( \bigO{ \sum_{\ell=1}^{d} \left( N_\text{rows}M m_\ell + m_\ell k_\ell^2 \right)} \) \\
\textbf{SeqSelect} &
\( \displaystyle \sum_{\ell=1}^{d} \cost_\text{cssp} \left(m_\ell, N_\text{rows} \cdot \frac{M}{m_\ell} \prod_{i=1}^{\ell-1} \frac{k_i}{m_i}, k_\ell \right) \) &
\( \bigO{ \sum_{\ell=1}^{d} \left( N_\text{rows}M m_\ell \prod_{i=1}^{\ell-1} \frac{k_i}{m_i} + m_\ell k_\ell^2 \right)} \) \\
\textbf{IterSelect} &
\( \displaystyle n_\text{iter} \sum_{\ell=1}^{d} \cost_\text{cssp} \left(m_\ell, N_\text{rows} \cdot \frac{K}{k_\ell}, k_\ell \right) \) &
\( \bigO{ n_\text{iter} \sum_{\ell=1}^{d} \left( N_\text{rows}K \frac{m_\ell^2}{k_\ell} + m_\ell k_\ell^2 \right)} \) \\
 &
\(\quad + (n_\text{iter}+1)\cdot \cost_\text{EIG}(K,N_\text{rows}) \) &
\(\quad + (n_\text{iter}+1)\cdot \cost_\text{EIG}(K,N_\text{rows}) \) \\
\bottomrule
\end{tabular}
\caption{Computational cost comparison for three structured CSSP methods. 
Each expression uses a general form where \( N_\text{rows} = N \) for deterministic CSSP and \( N_\text{rows} = K + p \) for the Sketch-First approach. 
While the second column gives the leading-order computational cost for the GKS-based implementation.}
\label{tab:cost_summary_vertical}
\end{table}

\paragraph{IterSelect}  
As shown in \eqref{def:unfoldIter}, for each iteration (\(1 \leq t \leq n_\text{iter}\)), the transpose of the \( \ell \)th unfolding of \( \T{Y}^{(t,\ell)} \) has shape
\[
(k_1 \cdots  k_{\ell-1}) \cdot (k_{\ell+1} \cdots k_d) \cdot N \times m_\ell = N{K}{k^{-1}_\ell}\times m_\ell.
\]
Thus, the total cost of this method is
\[
(n_\text{iter}+1)\cdot \cost_\text{EIG}(K,N)+  
n_\text{iter}\cdot \sum_{\ell=1}^{d} \cost_\text{cssp} \left(m_\ell, N{K}{k^{-1}_\ell}, k_\ell \right) \> \text{flops},
\]
where \( n_\text{iter} \) is the number of iterations needed for convergence.
\paragraph{Computational Cost of Evaluating the EIG}

The cost of evaluating the EIG objective in \eqref{def:Dopt} is denoted by
\[
\cost_\text{EIG}(n_c, n_r) = \mathcal{O}(n_c \cdot n_r \cdot \min\{n_c, n_r\}) \> \text{flops},
\]
which is dominated by the cost of computing the singular values of an \( n_r \times n_c \) matrix. The computational cost of the SVD is discussed in~\cite[Section 5.4]{demmel1997applied}, which we summarize here. 
When \( n_r > n_c \), the computation begins with a QR decomposition of the input matrix at a cost of \( \mathcal{O}(2n_c n_r - \tfrac{2}{3}n_c^3) \) flops, followed by the bidiagonalization of the resulting upper-triangular matrix \( \B{R} \), which requires \( \mathcal{O}(n_c^3) \) flops. The singular values of the bidiagonal matrix are then computed  with complexity \( \mathcal{O}(n^2_r ) \) flops.

We observe \( n_\text{iter} \leq 3 \) in all our numerical experiments.  
Since the cost of IterSelect depends on \( K \) rather than \( M \), this method can offer significant speedups compared to IndSelect and SeqSelect, when \( n_\text{iter} \) is small. 

\paragraph{Computational Cost of Sketch-First Approach}  
The cost follows directly from the deterministic case by replacing \( N \) with \( r=K + p \); see  \eqref{def:OmegaA}. The computational cost for both approaches is summarized in \Cref{tab:cost_summary_vertical}.

\section{Numerical Experiments} This section presents a comprehensive set of numerical experiments to evaluate the performance of the proposed method in the context of Bayesian inverse problems.
The aim is to demonstrate the versatility of the approach across a variety of relevant applications.
The experiments emphasize sensor placement strategies optimized for the EIG.

For each problem, we analyze the trade-offs between computational complexity and accuracy, benchmarking the proposed method against a greedy baseline and random designs. 
In addition to the deterministic version of our method, we  show the results for the sketch-based variant, Sketch-First, with an oversampling parameter $p=10$ for all experiments, namely, $r=K+10$. 
Specifically, we assess the effectiveness of the design variables selection procedure in a range of applications, including time-dependent partial differential equations (PDEs), seismic imaging, X-ray tomography, and flow reconstruction. 

All the numerical experiments were conducted on a Mac mini equipped with an Apple M1 processor, 8 GB of RAM, and 8 cores.

\subsection{Time-Dependent Problems}  
We consider a time-dependent Bayesian inference problem where the goal is to estimate the initial state of a dynamical system using observations from \( n_\text{s} \) fixed spatial sensors and a known evolution model \cite{Melina}.  

In this experiment, the underlying physical process is modeled by using the one-dimensional heat equation, which describes the evolution of temperature over space and time within a domain \( \Omega := (0,1) \). 
Although the true physical dynamic may be more complex, we assume that it can be adequately approximated by this model. 
The mathematical model governing the temperature distribution \( u(x,t) \) is given by
\begin{align}\label{eqn:1dheat}
\begin{aligned}
\frac{\partial u}{\partial t} &= 
\frac{\partial}{\partial x}  \left({\kappa} \frac{\partial u}{\partial x}   \right), & \text{in } \Omega \times (0, T), \\
u &= 0, & \text{on } \partial \Omega \times (0, T), \\
u &= u_0, & \text{on } \Omega \times \{0\}.
\end{aligned}
\end{align}
Here, $u(x, t)$ denotes the temperature at spatial position $x$ and time $t$, \( u_0 \) is the unknown initial temperature distribution, and ${\kappa}\equiv \sqrt{3}$ represents the thermal diffusivity. 
The inverse problem involves estimating the initial condition $u_0$ from discrete measurements of the state $u(x,t)$ in space and time at selected sensor locations. 
Observations are collected from a grid of \( n_s = 28 \) candidate sensors, uniformly distributed in \( \Omega \) away from the boundary, and corrupted with \( 2\% \) additive Gaussian noise to simulate measurement error.

\paragraph{PDE Discretization and Prior Covariance Matrix:}
We discretize the problem using the finite element method (FEM) with a Lagrange basis in space and the implicit Euler scheme in time. 
The spatial domain is resolved with \( 401 \) degrees of freedom, and \( 10 \) temporal snapshots are recorded at uniform intervals of \( 4 \times 10^{-3} \).
The mass and stiffness matrices  denoted by \( \B{N} \) and \( \B{K} \), respectively, define the prior covariance:  
\begin{align}
    \B{\Gamma}_\text{pr} = \left(\gamma\B{K} + \B{N}\right)^{-1} \B{N} \left(\gamma\B{K} + \B{N}\right)^{-1},
\end{align}  
where we set \( \gamma = 10^{-1} \).  

\paragraph{Matrices:} As described in \cite[Section 2]{Melina}, the system state at discrete time steps is given by  $
\B{u} = [\B{u}\t_0, \ldots, \B{u}\t_J]\t$,
where \( \B{u}_\ell \in \mathbb{R}^{n} \) represents the state at \( t = t_\ell \). The state evolves according to $\B{u}_{\ell} = \B{M}_{0,\ell}(\B{u}_0).$ The observation operator  \( \B{H} : \R^{n(J+1)} \to \R^{n_\text{s}(J+1)} \) acts independently at each time step and is given by    
\[
\B{H}(\B{u})= [\B{H}_0(\B{u}_0)\t, \ldots, \B{H}_0(\B{u}_J)\t]\t,
\]  
since the observation operator does not change in time. 
Here
\( \B{H}_0 : \R^n \to \R^{n_\text{s}} \) maps the state to observations at the \( n_\text{s} \) candidate sensor locations.
For simplicity we assume that the observation errors are uncorrelated across time, leading to the block-diagonal covariance matrix    
\begin{align}
    \B{\Gamma}_{\text{noise}} = \mathrm{blockdiag}(\B{R}_0, \ldots, \B{R}_N) \in \mathbb{R}^{n_\text{s}(J+1) \times n_\text{s}(J+1)},
\end{align}  
where \( \B{R}_\ell=\sigma^2_\ell\B{I} \in \mathbb{R}^{n_\text{s} \times n_\text{s}} \) represents the error covariance for observations at time \( t = t_\ell \), as described in~\cite[Section 6.1]{alexanderian2025optimal}. 
In this setting, the matrix $\B{A}$ takes the form
where  
\begin{align}
     \B{A}=  [\B{A}_0~ \cdots~ \B{A}_J], \quad \text{with } \B{A}_\ell = \sigma_{\ell}^{-1}\B{\Gamma}_\text{pr}^{\frac{1}{2}} \B{M}_{0,\ell}\t \B{H}_\ell\t.
\end{align}  

This results in the  matrix \( \B{A}\in \R^{401 \times (10\cdot 28)} \), which is then reshaped  into a three-dimensional tensor \( \T{X} \) of size \( 28 \times 10 \times 401 \), where  \( n_s=  28 \) corresponds to the total number of candidate sensor locations; \( T = 10 \) represents the number of discrete time snapshots; and  \( N = 401 \) is the dimension of the FEM space.

\paragraph{Experiment: EIG Structured Selection}
For this problem $M = n_sT$, and we want to perform subset selection of the form $\B{A}\B{S}=\T{X}_{(3)}\B{S} = \B{A}(\B{I}_T\otimes \B{S})$,
where the same \( k \) columns are selected from each block \( \B{A}_\ell \) for all \( \ell = 0, \dots, J \). This structure arises from the assumption that sensor locations are fixed in space. 
This is equivalent to setting \( \B{S}_2 = \B{I}_T \), or \( k_2 = T \), where the latter enforces a sorting of the snapshots by their relative contribution to EIG.
We consider two different values of $k$, namely, $5$ and $22$. For these values of $k$, we can determine the optimal sensor placement by an exhaustive search.

\begin{figure}[!hbt]
    \centering
    \begin{subfigure}[t]{0.45\linewidth}
        \centering
        \vspace{0pt}  
        \includegraphics[width=\linewidth]{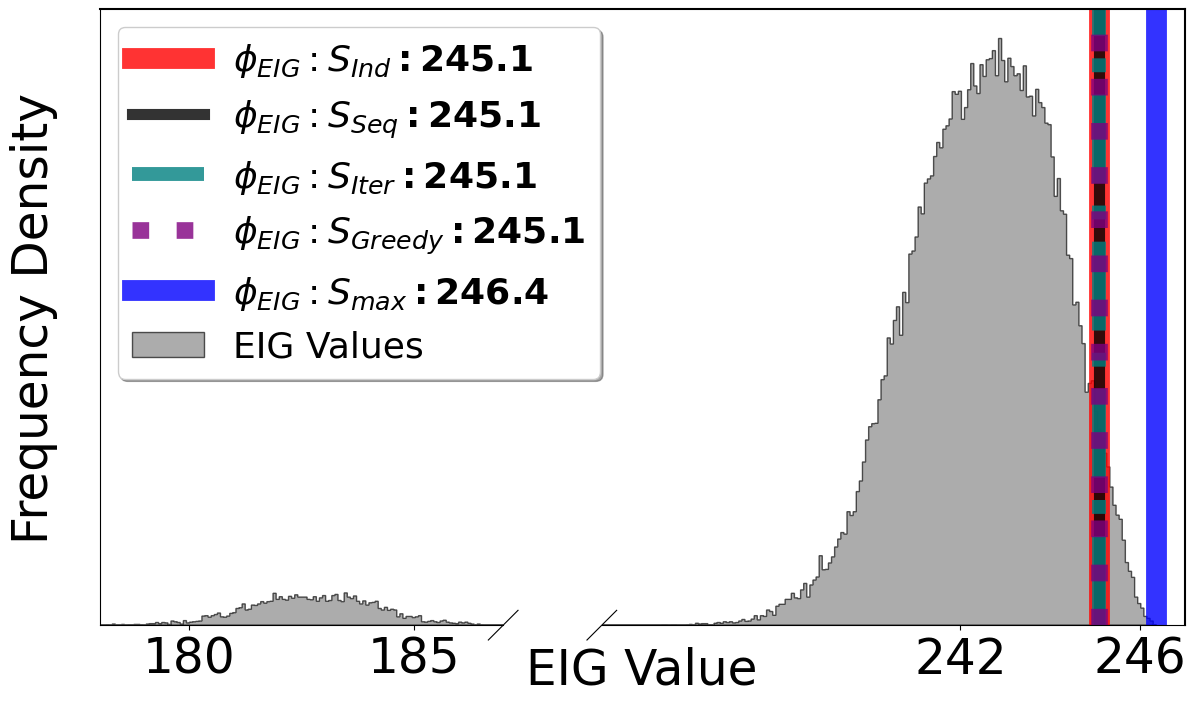}
        \caption{EIG values for selections with $k=5$.}
        \label{fig:SC4dVar5}
    \end{subfigure}
    \hfill
    \begin{subfigure}[t]{0.45\linewidth}
        \centering
        \vspace{0pt}  
        \includegraphics[width=\linewidth]{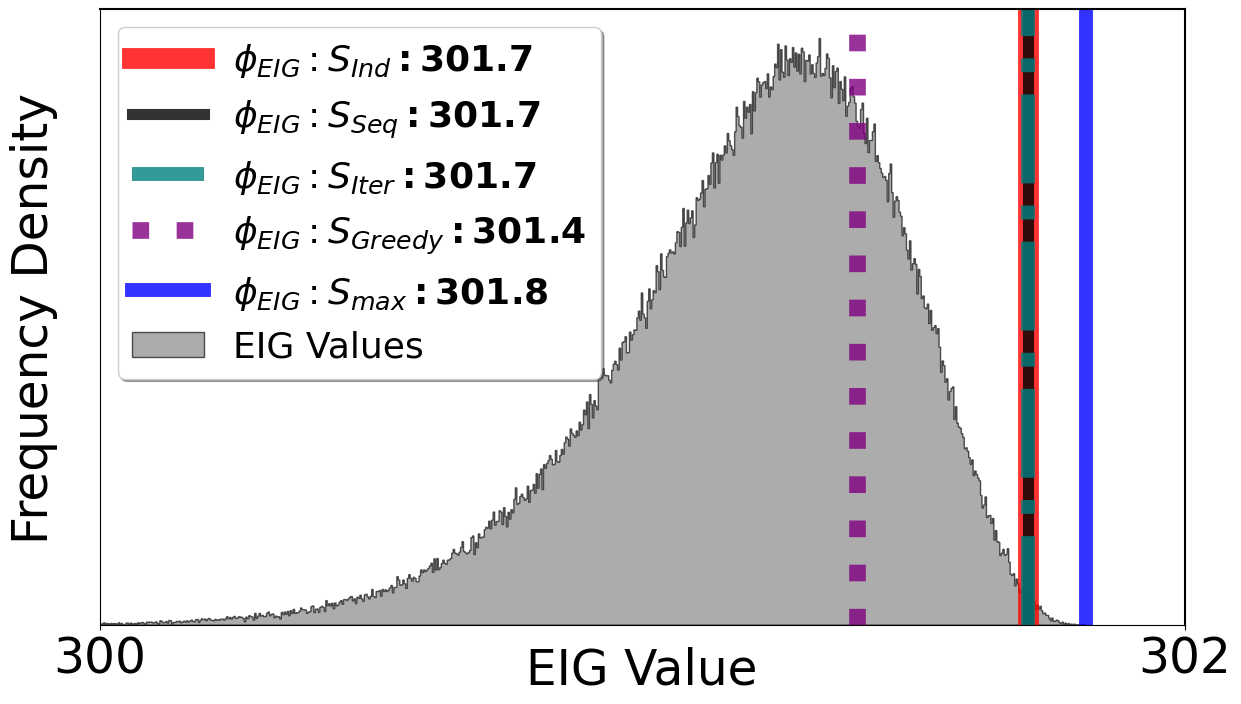}
        \caption{EIG values for selections with $k=22$.}
        \label{fig:SC4dVar22}
    \end{subfigure}
        \caption{Comparison of EIG values for different numbers of selected sensors ($k$) out of $n_s = 28$. All three proposed methods produce similar designs, achieving nearly identical EIG values.}

    \label{fig:SC4dVar}
\end{figure}

\paragraph{Results and Discussion}  
We compare the proposed methods against a greedy selection approach;  the results are displayed in Figure~\ref{fig:SC4dVar}. The histograms in the figure show the distributions of EIG values across all possible sensor configurations for both test cases, \( k = 5 \) and \( k = 22 \). Overlaid on each histogram are the EIG values achieved by the three structured selection methods: {IndSelect}, {SeqSelect}, and {IterSelect}.  

For \( k = 5 \), all three methods produce nearly identical sensor selections and achieve performance comparable to the best possible design.
Quantitatively, this performance is better than 96.13\% of all possible designs.
For comparison, we have also included the performance of the greedy selection method. 

For \( k = 22 \), all three proposed methods again perform comparably well and clearly outperform the greedy approach.
In fact, their performance exceeds that of 99.7\% of all possible designs.

We note the computational efficiency of {IterSelect}. For both cases  $k = 5$ and $k = 22$, {IterSelect} converged in only \( n_\text{iter} = 2 \) iterations.

\subsection{Seismic Tomography}\label{sec:seismic}
Seismic tomography is an imaging technique used to reconstruct the subsurface of the Earth by analyzing seismic waves. We use the IRTools toolbox~\cite{gazzola2019ir} to generate a synthetic problem from seismic tomography.  
In this problem instance the domain is \( (x,y) \in \Omega = [0,1]^2 \subset \mathbb{R}^2 \). 
We consider \( s = 32 \) sources, uniformly distributed along the right boundary (\( x = 1 \)), and \( q = 45 \) receivers, uniformly placed along the top boundary (\( y = 1 \)). Measurement noise is simulated by adding \( 2\% \) i.i.d. Gaussian noise to the data.
The design problem is to select \( k_1 \) sources (out of $s$) and \( k_2 \) receivers (out of $q$) that maximize the EIG. 
\paragraph{Discrete Problem and Prior Covariance Matrix:}  
We use the numerical implementation from~\cite{gazzola2019ir} and discretize the domain using a structured grid with \( n = 128 \times 128 \) nodes. The prior is derived from an elliptic problem and is given by \( \B\Gamma_{\rm pr}^{-1} = \gamma \B{K} \B{M}^{-1} \B{K} \), where \( \B{K} \) is the finite element representation of the discretized form of \( (\kappa^2 \text{id} - \Delta) \), where $\text{id}$ is the identity operator, with \( \gamma = \kappa^2 = 100 \) and \( \B{M} \) is the corresponding mass matrix.

\begin{figure}[ht!]
    \centering
    \begin{subfigure}{0.45\linewidth}
        \centering
        \includegraphics[width=\linewidth]{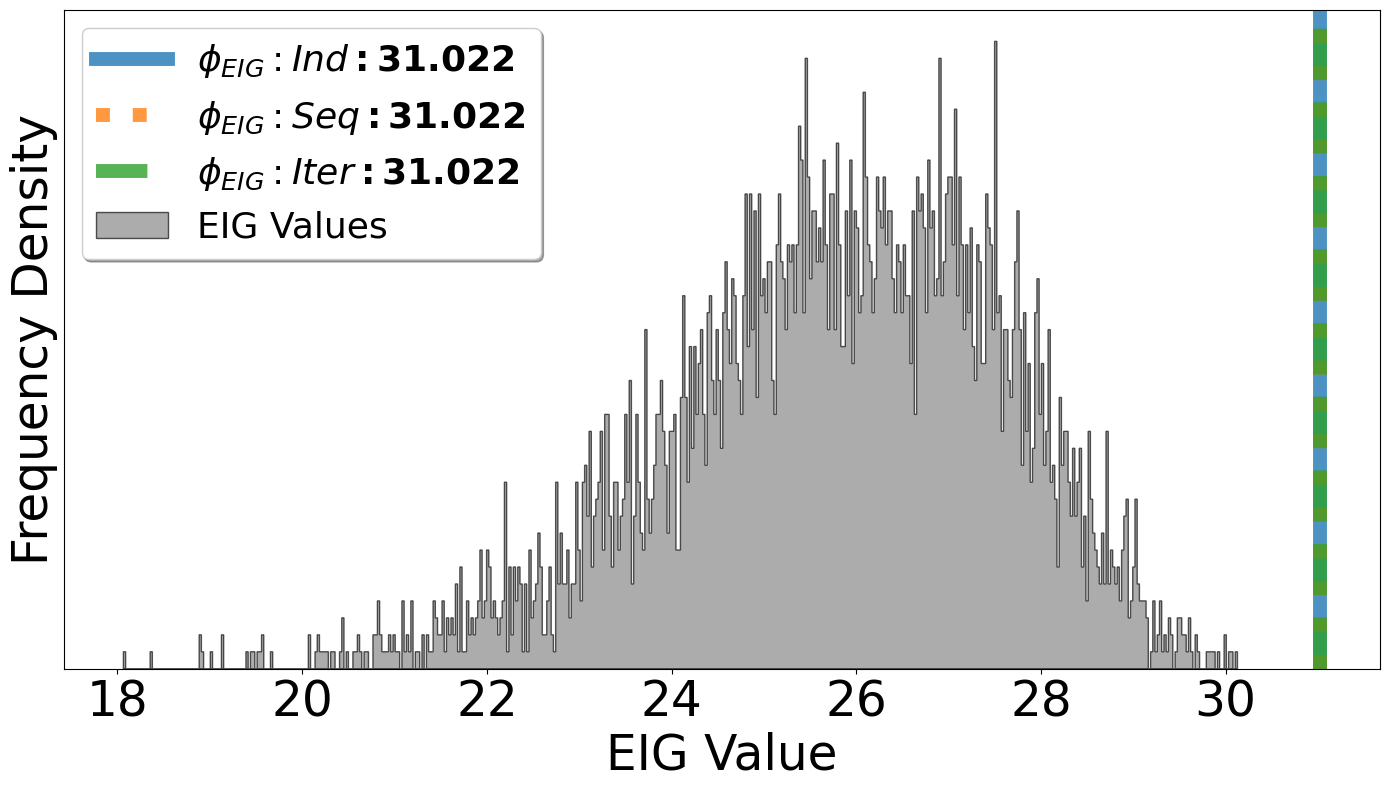}
        \caption{Comparison of EIG for the selection of \( k_1 = k_2 = 10 \).}
        \label{fig:dopt-seismic}
    \end{subfigure}
    \hfill
    \begin{subfigure}{0.45\linewidth}
        \centering
        \includegraphics[width=\linewidth]{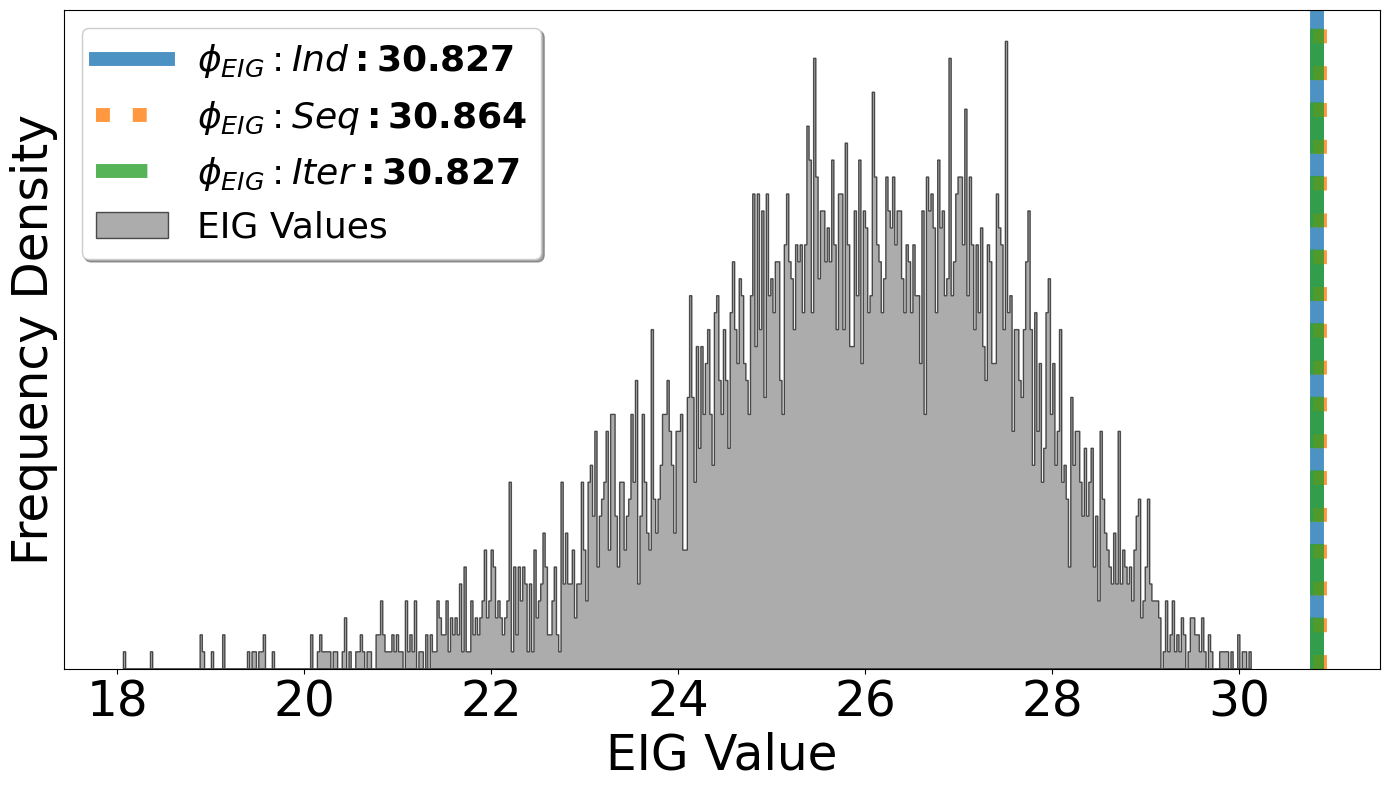}
        \caption{Comparison of EIG for sketching-based selection of \( k_1 = k_2 = 10 \).}
        \label{fig:map-seismic}
    \end{subfigure}
    \caption{Comparison of designs with and without sketching for \( k_1 = k_2 = 10 \) against  \( 5 \times 10^3 \) random source-receiver configurations sampled uniformly. 
    The Sketch-First method helps reduce computational costs while maintaining comparable performance. 
    }
    \label{fig:seismic-results}
\end{figure}

\paragraph{Experiment 1: EIG Structured Selection}  
In this setting the tensor \( \T{X} \in  \mathbb{R}^{32 \times 45 \times 4096} \) encodes the full set of sources, receivers, and grid points, respectively.  
The objective is to select \( k_1 \) sources and \( k_2 \) receivers, resulting in a structured design of the form \( \B{S} = \B{S}_2 \times \B{S}_1 \), where \( \B{S}_1 \) and \( \B{S}_2 \) represent selections of \( k_1 \) sources and \( k_2 \) receivers, respectively.  

We begin by evaluating {IndSelect}, {SeqSelect}, and {IterSelect} using \( k_1 = k_2 = 10 \), with the GKS method as the underlying CSSP algorithm.  
To benchmark their performance, we generate \( 5 \times 10^3 \) random source-receiver configurations sampled uniformly.

\paragraph{Results and Discussion}  

The proposed methods consistently outperform random designs and achieve higher EIG values.  
The left panel of Figure~\ref{fig:seismic-results} shows  EIG values for the sensor configurations selected by IndSelect, SeqSelect, and IterSelect.  
The right panel shows the results from the Sketch-First method, which works with a compressed tensor \( \T{Y} \in \mathbb{R}^{32 \times 45 \times 110} \).  
The IterSelect method required only two and three iterations, respectively, to reach convergence.

These results confirm that structured selection strategies are effective for this problem.  
They also show that sketching reduces computation while preserving accuracy in selecting informative source-receiver pairs.

\begin{figure}[!ht]
    \centering
    \begin{subfigure}{0.45\linewidth}
        \centering
        \includegraphics[width=\linewidth]{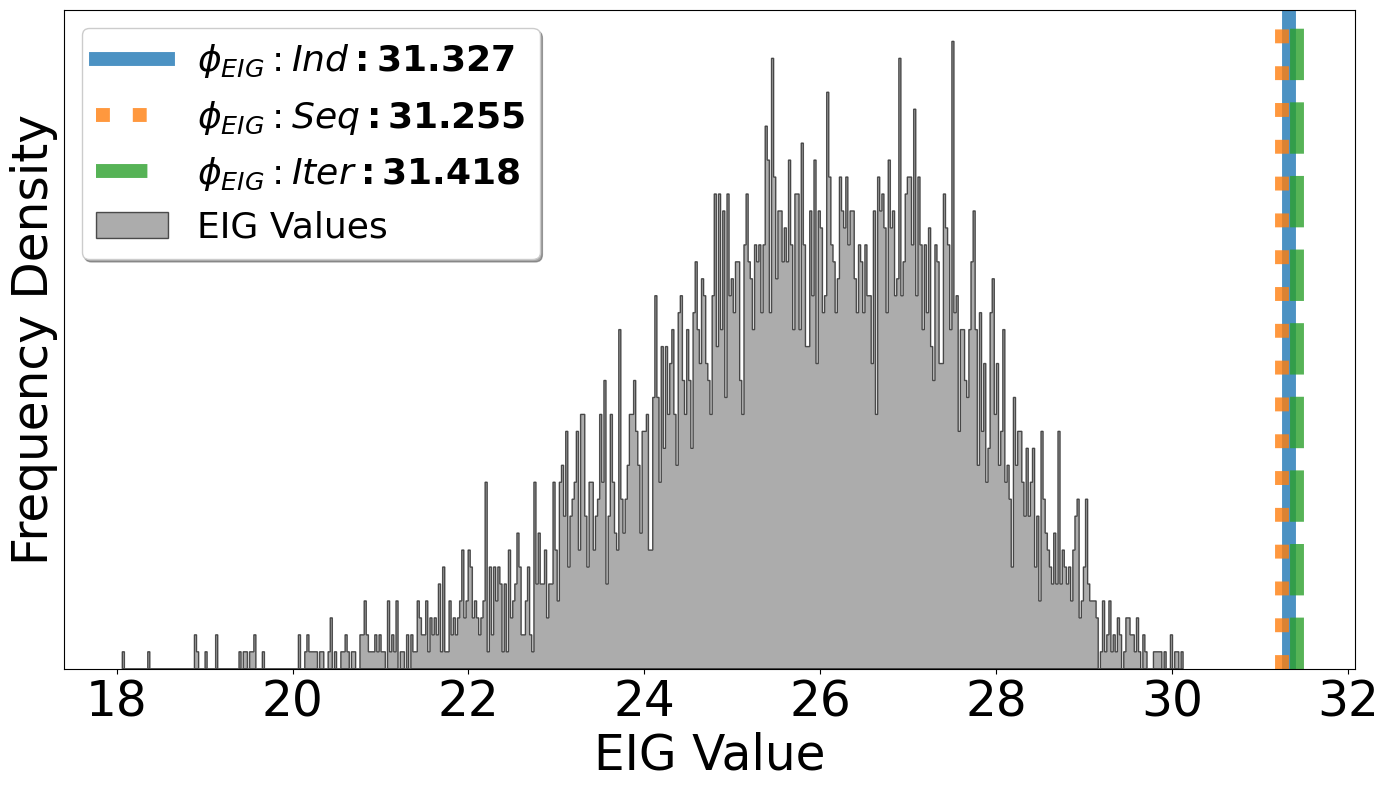}
        \caption{Greedy: \( k_1 = k_2 = 10 \) sensors.}
        \label{fig:dopt-seismic-greedy}
    \end{subfigure}
    \hfill
    \begin{subfigure}{0.45\linewidth}
        \centering
        \includegraphics[width=\linewidth]{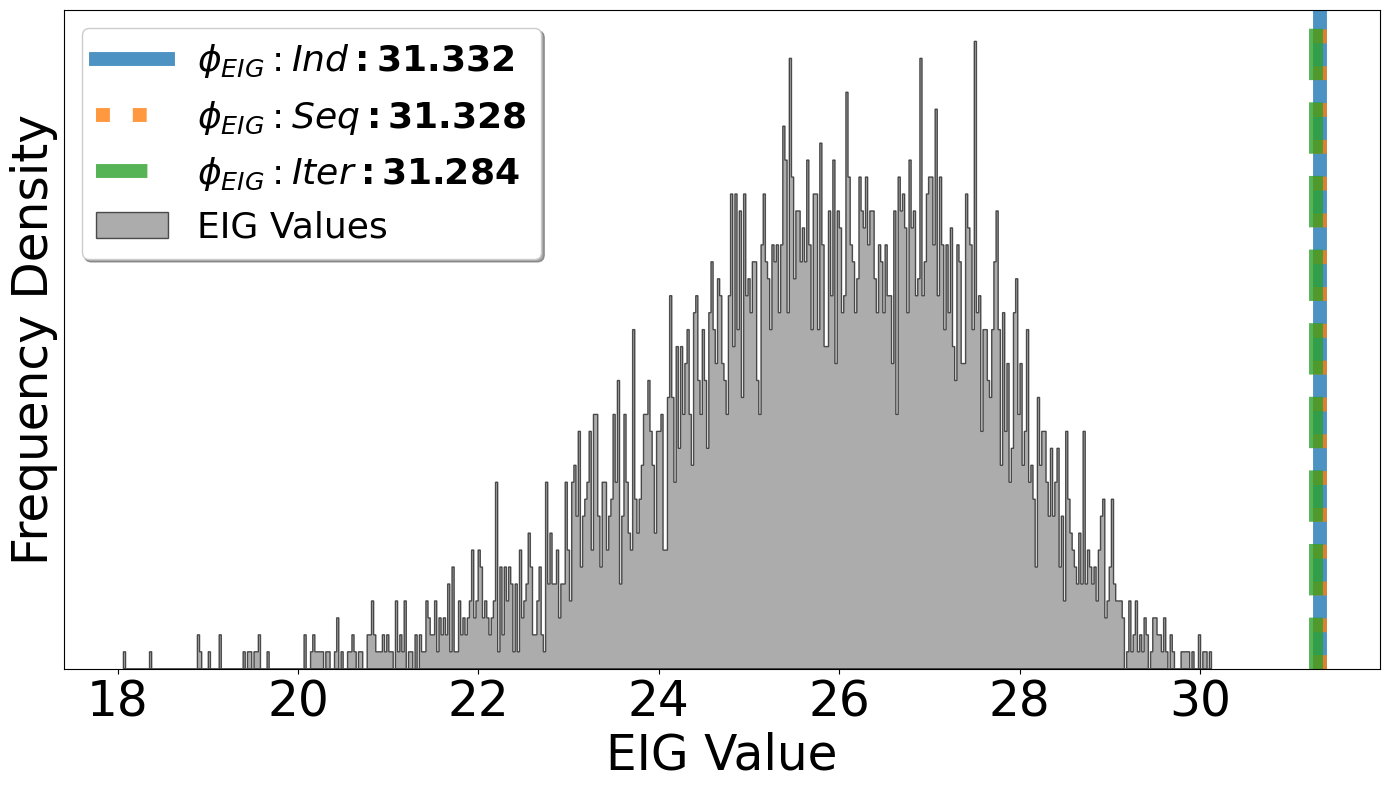}
        \caption{Greedy with sketching: \( k_1 = k_2 = 10 \).}
        \label{fig:map-seismic-greedy}
    \end{subfigure}
    \caption{Comparison of EIG values for sensor designs using the greedy method, with and without sketching, for \( k_1 = k_2 = 10 \).}
    \label{fig:seismic-results-greedy}
\end{figure}

\paragraph{Experiment 2: Greedy EIG Structured Selection}  
We now investigate greedy versions of {IndSelect}, {SeqSelect}, and {IterSelect} by replacing the CSSP step with a greedy search, as described in \Cref{sec:GreedySSSP}.  
We evaluate both the deterministic and Sketch-First variants of these methods under the same conditions as in the previous experiment.  
The goal is to assess the impact of the greedy heuristic on EIG performance while preserving the structured nature of the sensor selection.

\paragraph{Results and Discussion}  
\begin{table}[ht]
\centering
\renewcommand{\arraystretch}{1.2}
\begin{tabular}{lcccc}
\hline
\textbf{Method} & \textbf{GKS} & \textbf{GKS (Sketch-First)} & \textbf{Greedy} & \textbf{Greedy (Sketch-First)} \\
\hline
\textbf{IndSelect}  & \( 8.6349 \times 10^{-1} \) & \( 1.4253 \times 10^{-2} \) & 18.4519 & 0.5868 \\
\textbf{SeqSelect}  & \( 4.3040 \times 10^{-1} \) & \( 6.6851 \times 10^{-3} \) & 11.7843 & 0.3811 \\
\textbf{IterSelect} & \( 3.1077 \times 10^{-1} \) & \( 1.0841 \times 10^{-2} \) & 8.2991  & 0.3826 \\
\hline
\end{tabular}
\caption{Runtimes (in seconds) for selection methods using GKS and Greedy approaches, with and without sketching.}
\label{tab:runtime_horizontal}
\end{table}

Figure~\ref{fig:seismic-results-greedy} shows results using greedy selection instead of the GKS method, as in Figure~\ref{fig:seismic-results}.  
In both the deterministic and Sketch-First settings, greedy selection consistently outperforms all \( 5 \times 10^3 \) random designs and achieves slightly higher EIG values than the GKS-based methods do.  
The Sketch-First variants perform marginally better than their deterministic counterparts.  
The IterSelect method required three and four iterations, respectively, to converge.

Overall, IndSelect, SeqSelect, and IterSelect continue to perform well, matching the effectiveness of greedy strategies while remaining computationally efficient.  
The Sketch-First approach also proves effective in this setting. Across all methods, applying sketching consistently reduces runtimes by one to two orders of magnitude. 

In terms of runtime, the GKS-based methods are substantially faster compared to  their greedy counterparts, see Table~\ref{tab:runtime_horizontal}. 
Note that the greedy implementations were not optimized for speed; on the other hand, the GKS-based approaches use highly optimized linear algebra routines. While some care should be exercised while interpreting the timing results, this underscores a major advantage of the GKS-based approaches.

\paragraph{Experiment 3: Impact on Reconstruction Accuracy}

Next, we evaluate the reconstruction accuracy when restricting the parameter-to-observable map \( \B{F} \) and data \( \B{d} \) to the selected indices, as explained in \Cref{sec:SelectionMatrices} . We compare different structured selection methods against the full-data case (\( \B{S} = \B{I}_{45}\times \B{I}_{32} \)), where no selection is performed, and against \( 5 \times 10^3 \) random source-receiver configurations. 

\begin{figure}[ht!]
     \centering
     \includegraphics[width=0.6\linewidth]{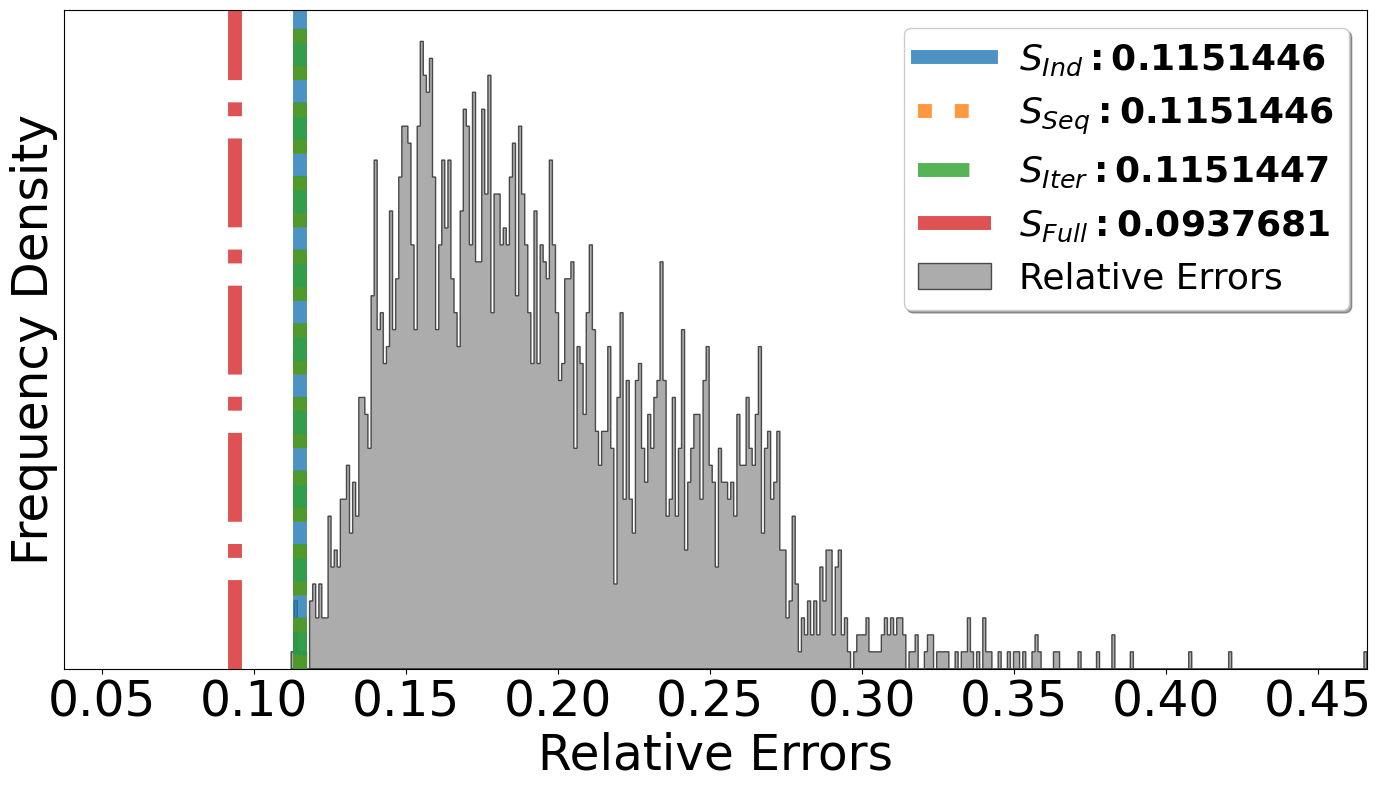}
     \caption{Relative reconstruction error for structured selection with \( k_1 = k_2 = 10 \).} \label{fig:rel-error-seismic}
\end{figure}

\paragraph{Results and Discussion}  
Figure~\ref{fig:rel-error-seismic} shows that using only 10 of the 32 sources and 10 of the 45 receivers—about \( 7\% \) of all possible source-receiver pairs—increases the relative error from \( 9.4\% \) to \( 11.5\% \).  
This small increase indicates that structured selection achieves a good balance between accuracy and computational savings.

\subsection{X-ray Tomography}  

X-ray tomography is a widely used imaging technique in which internal structures are inferred from a collection of X-ray projections taken at different angles. In a parallel-beam setup, each X-ray source emits rays that are detected across an array of evenly spaced detectors, resulting in line integrals of the attenuation field governed by the Radon transform \cite{Radon-1986}. In our experiment we consider 30 projection angles and 71 detector positions,  where each slice corresponds to a different projection angle and contains linear measurements of a \( 50 \times 50 \) attenuation map.

\paragraph{Discrete Problem and Prior Covariance Matrix:}  
The imaging domain is discretized into a regular grid of \( 50 \times 50 \) pixels, yielding 2,500 unknowns. The forward model is based on the discrete Radon transform, where each measurement corresponds to a line integral through the domain. The forward operator \( \B{F} \in \mathbb{R}^{2130 \times 2500}\) (corresponding to 30 angles and 71 detector positions per angle) and the unknown vector $\B{u} \in \R^{2500}$ model the attenuation coefficients.  We take the same type of prior covariance \(\B{\Gamma}_\text{pr}\)  as in the seismic experiment (Section~\ref{sec:seismic}).

\begin{figure}[!ht]
    \centering
    \begin{subfigure}{0.45\linewidth}
        \centering
        \includegraphics[width=\linewidth]{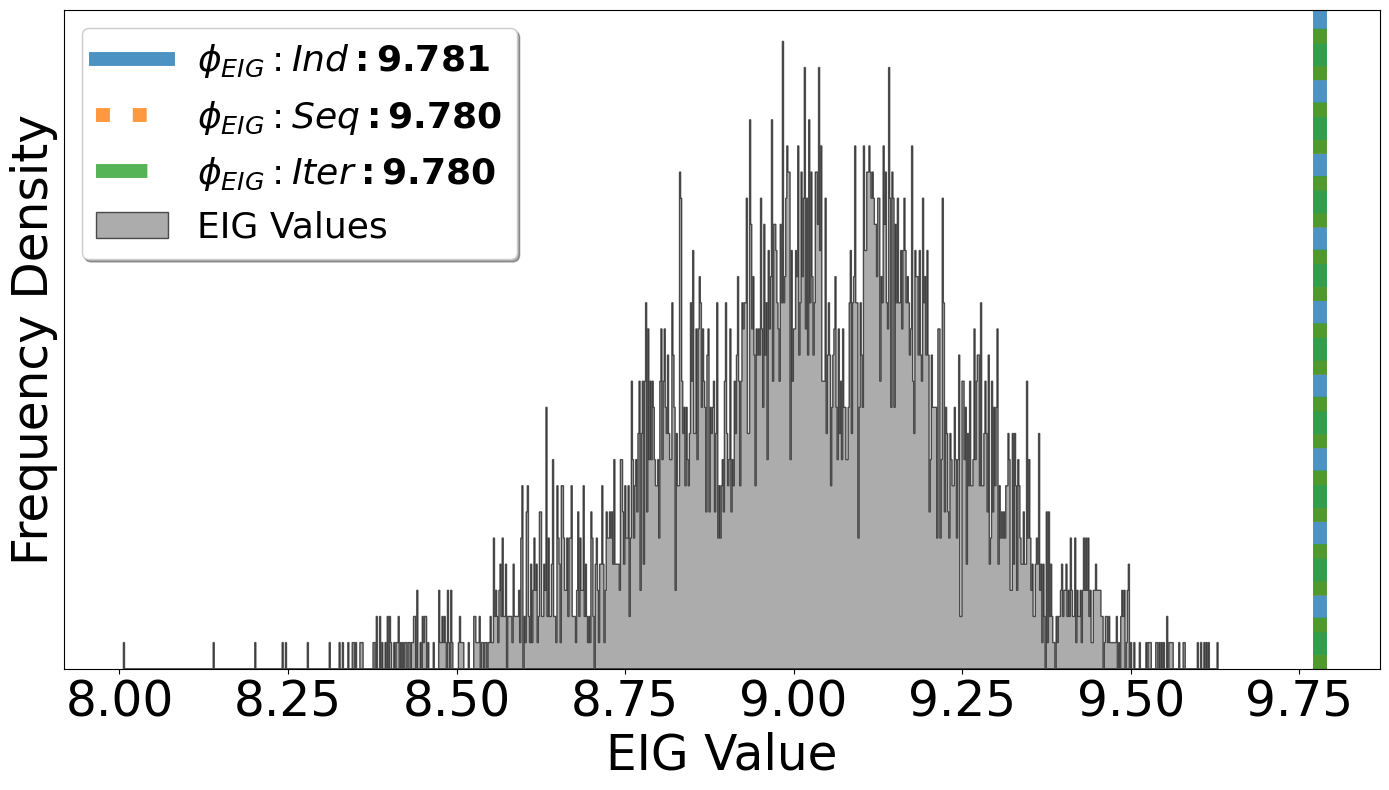}
        \caption{EIG values for different selection methods.}
        \label{fig:DoptXrays}
    \end{subfigure}
    \hfill
    \begin{subfigure}{0.45\linewidth}
        \centering
        \includegraphics[width=\linewidth]{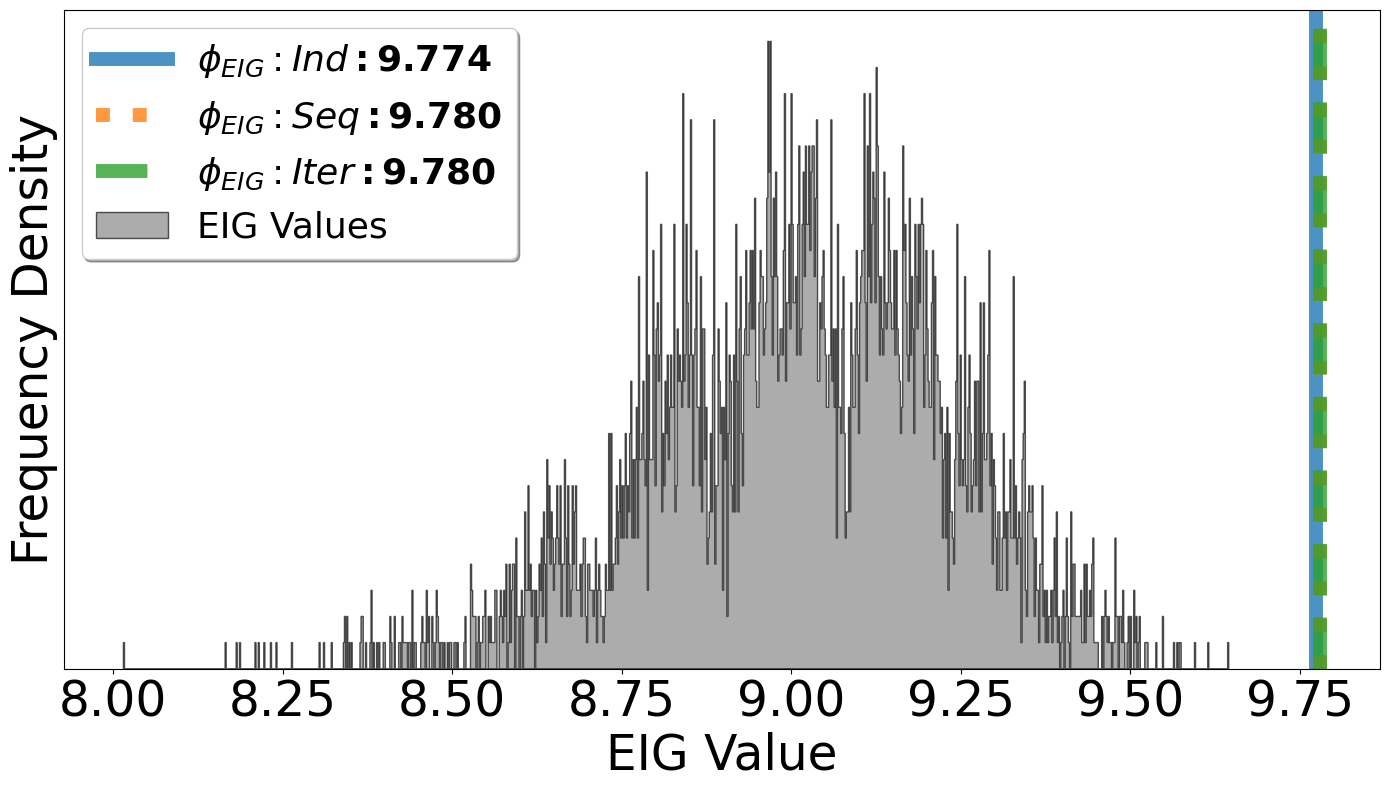}
        \caption{EIG values for the Sketch-First  method.}
        \label{fig:OtherImage}
    \end{subfigure}
    \caption{Comparison of EIG values in structured X-ray tomography design, with and without sketching.}
    \label{fig:DoptComparison}
\end{figure}

\paragraph{Experiment:}  
In this setting the corresponding  tensor  \( \T{X} \in \mathbb{R}^{30 \times 71 \times 2500} \)  encodes the full set of projection angles, detector positions, and grid points, respectively.
The design problem is to select the most informative \(k_1=20\) projection angles and \(k_2=25\) detector positions, as quantified by the EIG.
This yields a structured experimental design problem over the space of possible angle-detector pairs. 
\paragraph{Results and Discussion:}  
Figure~\ref{fig:DoptComparison} shows the results of different selection strategies. 
In panel~(a), we compare the proposed structured methods against the EIG values from $ 5.0 \times 10^3$ random designs. Panel~(b) shows the performance of the Sketch-First method, which deals with a tensor $\T{Y}\in \R^{30\times 71\times 510}$. 
IterSelect  required  two and three  iterations, respectively.

The proposed structured selection achieves EIG values better than any of the random designs.
Additionally, the sketching-based variant yields similar performance at significantly lower computational cost, demonstrating its effectiveness for large-scale problems.

\subsection{Flow Reconstruction}  

We consider a data-driven Bayesian inverse problem following the library-based framework in \cite{Lev}, applied to flow reconstruction tasks as in \cite{Tensor-based_flow}.  
Let \( \T{U} \in \mathbb{R}^{d_1 \times d_2 \times \cdots \times d_L} \) denote an unknown state vector discretized in a Cartesian grid in $L$ dimensions.  We assume that \( \T{U} \) can be well approximated on a reduced basis as  
\[
\operatorname{vec}(\T{U}) \approx \B{\Phi m},
\]
where \( \B{\Phi} \in \mathbb{R}^{M \times N} \) is a modal basis, \( \B{m} \in \mathbb{R}^N \) are the corresponding reduced coordinates, and $M=d_1\cdots d_L$.  
We consider noisy (vectorized) observations \( \B{y} \in \mathbb{R}^M \), assumed to be related to the reduced state via the model  
\begin{align}\label{def:dataFlowRec}
    \B{y} = \B{\Phi m} + \B{\eta},
\end{align}
where \( \B{\eta} \sim \mathcal{N}(\B{0}, \sigma_{\B{R}}^{2}\B{I}) \). The inverse problem involves recovering the coefficients $\B{m}$ from the data $\B{y}$. The rows of $\B\Phi$ map to the candidate sensor locations; we can place a limited number of sensors, which we also place in a Cartesian grid. The OED problem involves finding the optimal sensor locations. 

\paragraph{Dataset and Preprocessing}  
We use the Tangaroa flow dataset, which captures turbulent airflow around a NIWA research vessel \cite{Popinet04:Tangaroa}. As in~\cite{Tensor-based_flow}, we only focus on recovering the velocity component in the $x$-direction. The data is organized as a fourth-order tensor. The first three dimensions represent spatial coordinates (\(x, y, z\)), and the fourth dimension corresponds to time snapshots.  The resulting tensor \( \T{X} \) has shape \( 300 \times 180 \times 120 \times 201 \).

To reduce computational cost, we downsample the spatial dimensions by a factor of 2. We then center the dataset by subtracting the temporal mean at each spatial location. Specifically, we compute the mean tensor \( \overline{\T{X}}^\text{tr} \in \mathbb{R}^{150 \times 90 \times 60} \) over the first $t=150$ snapshots and subtract it from each corresponding frame. The result is a centered training tensor \( \T{X}^\text{tr} \in \mathbb{R}^{150 \times 90 \times 60 \times 150} \) and a test tensor \( \T{X}_\text{test} \in \mathbb{R}^{150 \times 90 \times 60 \times 51} \) containing the remaining 51 time steps. We add $2\%$ of Gaussian noise to the data to represent measurement noise. 

\paragraph{Basis and Prior:} In this setting, we construct the basis \( \B{\Phi}= \) from the
left singular vectors of the centered snapshot matrix
$ \B{X}_\text{snap} = (\T{X}_{(4)}^\text{tr})\t \in\R^{ 150\times (150\cdot 90\cdot 60)}$, whose SVD is 
denoted \( \B{X}_\text{snap} = \B{U\Sigma V}\t\). The mean snapshot is subtracted prior to decomposition. The authors in \cite{Lev} propose to consider as a prior covariance for $\B{m}$ given by
\begin{align}
   \B{\Gamma}_\text{pr}= \frac{1}{t-1}\B{\Sigma}^2,
\end{align}
where $t$ is the number of snapshots in the training tensor. We adopt this prior covariance in our numerical experiments.

\begin{figure}[ht!]
    \centering
    \includegraphics[width=0.6\linewidth]{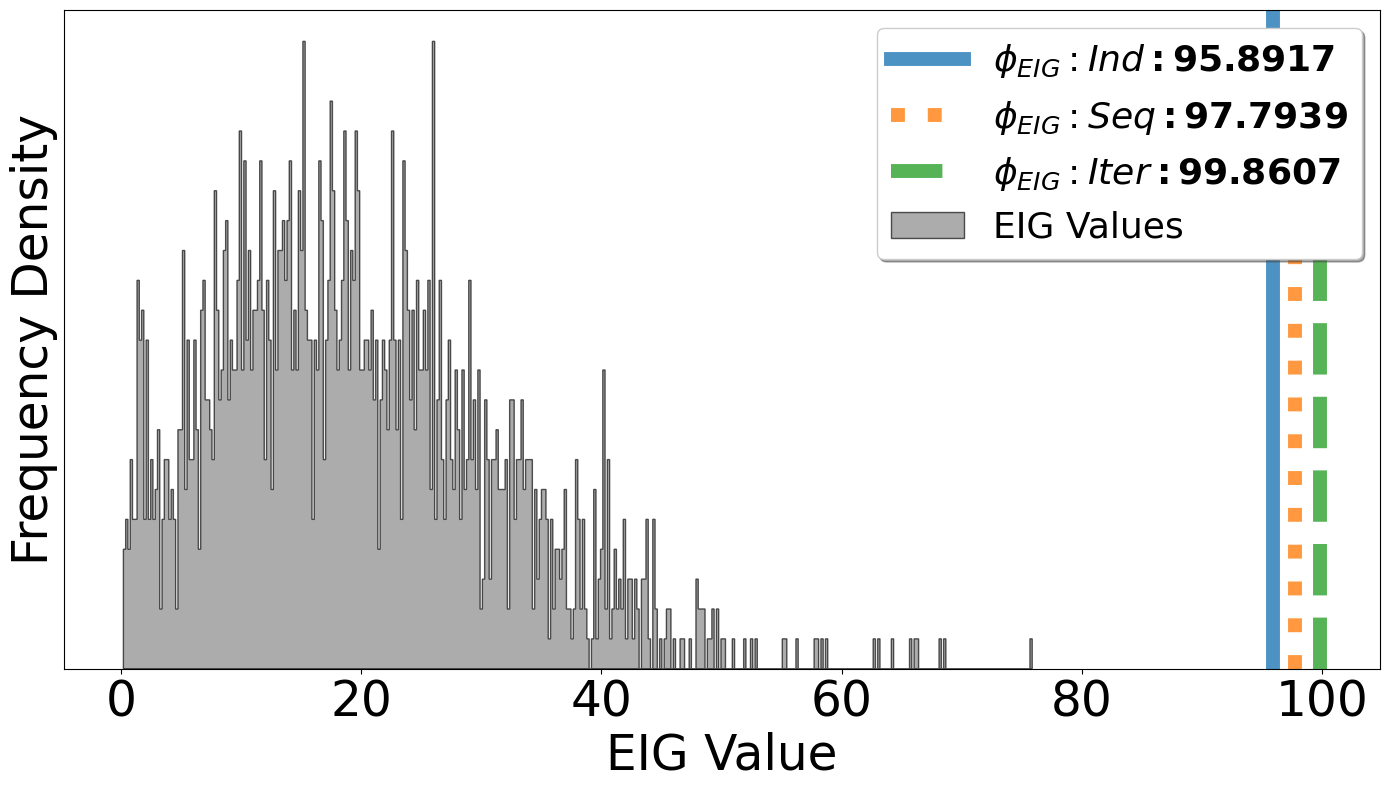}
    \caption{Comparison of EIG values of structured CSSP selection strategies on the flow reconstruction application. }
    \label{fig:DoptDEIM}
\end{figure}

\paragraph{Results and discussion:}

In this experiment, we assess the effectiveness of different selection strategies for identifying informative variables along each spatial mode of the Tangaroa dataset. Specifically, we apply the \textit{IndSelect}, \textit{SeqSelect}, and \textit{IterSelect} methods to select  \( k_1 = k_2 = k_3 = 5 \), giving a total of $5^3 = 125$ sensors. For comparison, we also evaluate the EIG values obtained from $2 \times 10^3 $ random designs of the same size. The results are summarized in Figure~\ref{fig:DoptDEIM}.

\begin{table}[!ht]
\centering
\renewcommand{\arraystretch}{1.2}
\begin{tabular}{lcc}
\toprule
\textbf{Method} & \textbf{GKS} & \textbf{Greedy} \\
\midrule
\textbf{IndSelect} & 68.8012 & 321.4451 \\
\textbf{SeqSelect} & 32.6082 & 81.9791 \\
\textbf{IterSelect} & 0.5883 & 4.7337 \\
\bottomrule
\end{tabular}
\captionsetup{width=\textwidth}
\caption{Runtimes (in seconds) for GKS-based and Greedy-based selection methods. }
\label{tab:runtime_comparison_flow}
\end{table}
Our selection strategies consistently outperform the random designs. Notably, \textit{IterSelect} achieves the best performance in terms of EIG.

Additionally, as shown in Table~\ref{tab:runtime_comparison_flow}, the GKS-based implementations are significantly more efficient computationally than their greedy-based counterparts. 
Among the methods, IterSelect is by far the fastest. This is largely due to the fact that, in essence, IterSelect operates on a small tensor, since the number of target sensors is only $5^3$ in this experiment.

\section{Conclusions and Discussion}
In this article, we developed new algorithms for optimal experimental design via structured column subset selection. This generalizes the CSSP approach for $d = 1$, for structured designs. The proposed methods used tensor decompositions and CSSP applied to different mode unfoldings.  Several methods are proposed in this paper, which have partial analogs in terms of tensor decompositions.  The methods also leveraged randomized techniques to efficiently handle structured column subset selection. In numerical experiments, the proposed methods produce near-optimal EIG designs. The proposed methods are more computationally efficient, with up to $50\times$ speed ups compared to greedy approaches, and give comparable or better designs. The range of applications extends the  OED using CSSP in~\cite{eswar2023optimal} to a much wider class of problems, thereby broadening the impact of the work. 

There are several avenues for future work. Efficient computation of tensor decompositions is an active area of research. Using the connection to tensor decompositions established in this paper, it may be possible to derive new algorithms for structured column selection. Similarly, we only used a specific form of randomization in this paper; however, many other variants are possible, which we can explore. From a theoretical perspective, we were not able to derive analogs of Lemma~\ref{lemma:csspoed} for the structured column selection. Future work could explore the analysis of these algorithms. Finally, extensions to nonlinear OED problems and for other criteria are also of interest and are currently under investigation.

\bmhead{Supplementary information} Not applicable.
\bmhead{Acknowledgements}
HD and AKS were supported by the Department of Energy, Office of Science,
Advanced Scientific Computing Research (ASCR) Program through the awards DE-SC0023188 and DE-SC0025262.
SE, VR, and ZWD were supported
by the U.S. Department of Energy, Office of Science, ASCR Program under contracts DE-AC0206CH11357 and DE-SC0023188.
\bmhead{Competing Interests} The authors have no competing interests to declare.

\bibliography{refs}

\begin{appendices}
\section{Deterministic Structured Greedy Algorithms}\label{sec:GreedySSSP}

The selection strategies {IndSelect}, {SeqSelect}, and {IterSelect} are versatile and can be used with any CSSP method, not limited to the GKS method; see \Cref{sec:SSSPAlg} for the general framework.

In this section we present a structured greedy alternative that aligns with the CSSP interface. The pseudocode below outlines a greedy selection algorithm that can be used as a drop-in replacement for the  CSSP subroutine.
While a greedy selection can be applied directly to the matrix $\B{A}$ defined in \eqref{def:A} or its sketched version $\B{\Omega A}$ in \eqref{def:OmegaA}, our structured variant is tailored for tensor-based selection.
\Cref{alg:GreedyIndselect} presents a greedy CSSP algorithm for the IndSelect method.

\begin{algorithm}[!ht]
\caption{\textbf{IndSelect: Independent Mode Selection via Greedy CSSP}}
\label{alg:GreedyIndselect}
\begin{algorithmic}[1]
\REQUIRE Matrix \( \B{A} \in \R^{N \times M} \), number of selected indices \( \B{k} = (k_1, \dots, k_d) \)
\ENSURE Selection matrices \( (\B{S}_1, \dots, \B{S}_d) \) for each mode

\STATE Reshape matrix \( \B{A} \) into a tensor \( \T{X} \in \R^{m_1 \times \cdots \times m_d \times N} \)

\FOR{\( \ell = 1,\dots,d \)}  
    \STATE \# Compute mode-\( \ell \) unfolding (transposed)  
    \STATE \( \B{B} \gets \B{X}\t_{(\ell)} \in \R^{(m_1 \cdots m_{\ell-1} \cdot m_{\ell+1} \cdots m_d \cdot N) \times m_\ell} \)
    
    \STATE Initialize \( \B{B}_S \gets \B{0}_{(NM m_\ell^{-1}) \times k_\ell} \), \( \mathcal{S} \gets \B{0}_{1 \times k_\ell} \)
    
    \FOR{\( j = 1,\dots,k_\ell \)}  
        \STATE \( \text{best\_gain} \gets -\infty \), \( \text{best\_index} \gets -1 \)
        \STATE \( \text{shape} \gets (m_1,\dots,m_{\ell-1}, j, m_{\ell+1},\dots,m_d, N) \)
        
        \FOR{\( i = 1,\dots,m_\ell \)}  
            \IF{\( i \in \mathcal{S} \)}  
                \STATE \textbf{continue}  
            \ENDIF  

            \STATE Copy column \( i \): \( \B{B}_S(:, j) \gets \B{B}(:, i) \)
            \STATE Reshape: \( \T{B}^{\text{temp}} \gets \text{reshape}(\B{B}_S(:,1:j), \text{shape}) \)
            \STATE Evaluate design: \( \Psi_{\text{temp}} \gets \Psi\left(\B{B}_{(d+1)}^{\text{temp}}\right) \)

            \IF{\( \Psi_{\text{temp}} > \text{best\_gain} \)}  
                \STATE \( \text{best\_gain} \gets \Psi_{\text{temp}} \), \( \text{best\_index} \gets i \)  
            \ENDIF  
        \ENDFOR  

        \STATE \( \mathcal{S}(j) \gets \text{best\_index} \), \( \B{B}_S(:, j) \gets \B{B}(:, \text{best\_index}) \)
    \ENDFOR  

    \STATE Set selection matrix: \( \B{S}_\ell \gets \B{I}(:, \mathcal{S}) \)
\ENDFOR  

\RETURN \( (\B{S}_1, \dots, \B{S}_d) \)
\end{algorithmic}
\end{algorithm}

\subsection{Derivation of Computational Cost}
 We now focus on computational cost the IndSelect approach; we start by analyzing the selection along the first mode. In the initial step, the algorithm evaluates the EIG criterion for each of the $m_1$ slices independently. Each slice has size $1 \times m_2 \times \cdots \times m_d \times N$, so the total cost of this step is
\[
m_1 \cdot \cost_\text{EIG}(1 \cdot m_2 \cdots m_d, N) \> \text{flops}.
\]
In the second step, we evaluate the EIG criterion for each of the remaining $m_1 - 1$ candidates.
Each evaluation now uses a subtensor with size $2 \times m_2 \times \cdots \times m_d \times N$.
This process continues, increasing the number of combined slices at each step, until $k_1$ slices are selected. 
At step $\rho$, the algorithm considers $m_1 - \rho + 1$ candidates, each of size $\rho \times m_2 \times \cdots \times m_d \times N$. Therefore, the total cost for selection along the first mode is
\[
\sum_{\rho=1}^{k_1} (m_1 - \rho + 1) \cdot \cost_\text{EIG}(r \cdot m_2 \cdots m_d, N) \> \text{flops}.
\]
Extending this to all $d$ modes and letting $M = m_1 \cdot m_2 \cdots m_d$, the total cost across all modes is approximately
\[
\sum_{\ell=1}^d \sum_{\rho=1}^{k_\ell} (m_\ell - \rho + 1) \cdot \cost_\text{EIG}\left(\rho \cdot \frac{M}{m_\ell}, N\right) \> \text{flops}.
\]

\( \)

\subsection{Summary of Computational Costs: Greedy Approaches}

The computational cost analysis for each greedy method follows similar reasoning to the discussion presented earlier in this section. Therefore, we omit the step-by-step derivations and summarize the total costs in \Cref{tab:cost_Greedy}. In this summary, we assume that the modes are processed in the order $1, 2, \ldots, d$.
To simplify the analysis of the greedy approach, we assume that
$
M \le N,
$ 
which ensures that for each \(1 \le \ell  \le d\) and all \(1 \le \rho  \le k_\ell\), the following cost estimate holds:
\[
\cost_{\mathrm{EIG}}\left(\rho \cdot \frac{M}{m_\ell}, N\right) = \mathcal{O}\left(N \left(\rho \cdot \frac{M}{m_\ell}\right)^2\right) \> \text{flops}.
\]

\begin{table}[ht]
\centering
\renewcommand{\arraystretch}{1.6}
\begin{tabular}{@{}lll@{}}
\toprule
\textbf{Method} & \textbf{Computational Cost} & \textbf{Explicit Cost} (leading-order terms) \\
\midrule
\textbf{IndSelect} &
\( \displaystyle \sum_{\ell=1}^d \sum_{\rho=1}^{k_\ell} (m_\ell - \rho + 1) \cdot \cost_\text{EIG}\left(\rho \cdot \frac{M}{m_\ell}, N\right) \) &
 $\displaystyle NM^2 \sum_{\ell=1}^d \sum_{\rho=1}^{k_\ell} (m_\ell - \rho + 1)\frac{\rho^2}{m_\ell^2}$ \\
\textbf{SeqSelect} &
\( \displaystyle \sum_{\ell=1}^d \sum_{\rho=1}^{k_\ell} (m_\ell - \rho + 1) \cdot \cost_\text{EIG}\left(\rho \cdot \frac{M}{m_\ell} \prod_{i=1}^{\ell-1} \frac{k_i}{m_i}, N\right) \) &
 $\displaystyle NM^2 \sum_{\ell=1}^d \sum_{\rho=1}^{k_\ell} (m_\ell - \rho + 1)\frac{\rho^2}{m_\ell^2} \prod_{i=1}^{\ell-1} \frac{k^2_i}{m^2_i}$ \\
\textbf{IterSelect} &
\( \displaystyle n_\text{iter} \sum_{\ell=1}^{d}  \sum_{\rho=1}^{k_\ell} (m_\ell - \rho + 1) \cost_\text{EIG}\left( \rho\cdot \frac{K}{k_\ell}, N \right) \) \newline
 &
$\displaystyle n_\text{iter}\cdot  NK^2 \sum_{\ell=1}^d \sum_{\rho=1}^{k_\ell} (m_\ell - \rho + 1)\frac{\rho^2}{k_\ell^2}$ \\
& $+ (n_\text{iter} +1) \cdot \cost_\text{EIG}(K,N)$ & $+ (n_\text{iter} +1)NK^2$ \\ 
\bottomrule
\end{tabular}
\caption{Computational cost of structured CSSP using greedy selection methods, assuming the input matrix $\B{A}$ is precomputed.
}
\label{tab:cost_Greedy}
\end{table}

\paragraph{Computational cost comparison between the greedy and GKS-based methods}

To simplify the analysis, we assume that: \(k_1 = \dots = k_d = k\) and \(m_1 = \dots = m_d = m\). 
This allows for a clearer comparison of the leading-order computational costs of the greedy and GKS-based methods. 
We summarize the leading-order costs in Table~\ref{tab:cost-gks-greedy}. 

\begin{table}[!ht]
\centering
\renewcommand{\arraystretch}{1.7}
\begin{tabular}{l|c|c}
\textbf{Method} & \textbf{GKS-based Cost} & \textbf{Greedy Cost} \\
\hline
\textbf{IndSelect} & 
  $d(Nm^{d+1} + m k^2)$ & 
  $\displaystyle d \cdot N m^{2d-2} \cdot \left( \frac{(m+1)k^3}{3} - \frac{k^4}{4}\right)$ \\
\textbf{SeqSelect} & 
  $\displaystyle  \frac{1 - \left(\frac{k}{m}\right)^{2d}}{1 - \frac{k^2}{m^2}} \cdot Nm^{d+1}   + dm k^2$ & 
  $\displaystyle \frac{1 - \left(\frac{k}{m}\right)^{2d}}{1 - \frac{k^2}{m^2}} \cdot N m^{2d-2} \cdot \left( \frac{(m+1)k^3}{3} - \frac{k^4}{4}\right)$ \\
\textbf{IterSelect} & 
  $\displaystyle n_\text{iter} \cdot d(Nk^{d-1} m^2 + m k^2)$ & 
  $\displaystyle n_\text{iter} \cdot d N k^{2d-2} \cdot \left( \frac{(m+1)k^3}{3} - \frac{k^4}{4} \right)$ \\
 & 
  $\displaystyle + (n_\text{iter}+1)Nk^{2d}$ & 
  $\displaystyle + (n_\text{iter}+1)Nk^{2d}$ 
\end{tabular}
\caption{Leading-order terms from computational cost for the GKS-based and greedy approaches across three methods, under the assumption that \(k_1 = \dots = k_d = k\) and \(m_1 = \dots = m_d = m\).}
\label{tab:cost-gks-greedy}
\end{table}

It is evident from the table that the GKS-based approaches yields significantly lower computational complexity than the corresponding greedy methods. The speedups are particularly high if $k \ll m$ and the dimension $d$ is high. 
\end{appendices}

\end{document}